\documentclass[12pt]{article}
\usepackage{framed}
\usepackage{amsmath}
\usepackage{amsthm}
\usepackage{amsfonts}
\usepackage[utf8]{inputenc}
\usepackage{color}
\usepackage[T1]{fontenc}
\usepackage[margin=2.5cm]{geometry}
\usepackage{stmaryrd}
\usepackage{graphicx}
\usepackage{dsfont}
\usepackage{mathrsfs}
\usepackage{combelow}
\usepackage{placeins}
\usepackage[font=small,labelfont=bf]{caption}

\renewcommand{\Re}{\operatorname{Re}}
\renewcommand{\Im}{\operatorname{Im}}

\theoremstyle{plain}

\usepackage[symbol]{footmisc}

\newtheorem{theo}{Theorem}[section]
\newtheorem{lemme}[theo]{Lemma}
\newtheorem{prop}[theo]{Proposition}

\renewcommand{\sp}{\hspace{0.2cm}}

\newcommand{\dd}{\,\mathrm{d}}
\newcommand{\R}{\mathbb{R}}
\newcommand{\C}{\mathbb{C}}

\newcommand{\ds}{\displaystyle}

\newcommand{\eps}{\varepsilon}
\newcommand{\supp}{\operatorname{supp}}

\newcommand{\Z}{\mathbb{Z}}
\newcommand{\vs}{\vspace{1mm}}
\newcommand{\der}[2]{\frac{\dd #1}{\dd #2}}

\newcommand{\tgm}{\widetilde{\gamma}}

\newcommand{\vertiii}[1]{{\left\vert\kern-0.25ex\left\vert\kern-0.25ex\left\vert #1 
    \right\vert\kern-0.25ex\right\vert\kern-0.25ex\right\vert}}

\numberwithin{equation}{section}

\title{Long time confinement of vorticity around a stable stationary point vortex in a bounded planar domain}
\author{Martin Donati and Dragoș Iftimie}
\date{}

\begin{document}
\maketitle
\begin{abstract}
    In this paper we consider the incompressible Euler equation in a simply-connected bounded planar domain. We study the confinement of the vorticity around a stationary point vortex. We show that the power law confinement around the center of the unit disk obtained in \cite{marchiorobutta2018} remains true in the case of a stationary point vortex in a simply-connected bounded domain. The domain and the stationary point vortex must satisfy a condition expressed in terms of the conformal mapping from the domain to the unit disk. Explicit examples are discussed at the end.
\end{abstract}

\section{Introduction and main result}

To study the behavior of an incompressible inviscid fluid, we consider the planar Euler equations in a bounded domain $\Omega \subset \R^2$:
\begin{equation}
    \begin{cases}\label{Euler} \partial_t u(x,t) + u(x,t)\cdot \nabla u(x,t) = - \nabla p(x,t), & \forall(x,t) \in \Omega\times\R^*_+ \vs \\  u(x,0) = u_0(x), & \forall x \in \Omega \\ \nabla \cdot u (x,t) = 0,& \forall(x,t) \in \Omega\times\R_+ \\ u(x,t)\cdot \vec{n} = 0,&\forall (x,t) \in \partial\Omega\times\R_+\end{cases}
\end{equation}
where $u$ denotes the velocity of the fluid, $p$ its internal pressure, and $n$ is the exterior normal to $\partial\Omega$. If $\Omega = \R^2$ the boundary condition should be changed into a vanishing condition at infinity. We define the fluid's vorticity by $\omega = \partial_1 u_2 - \partial_2 u_1$, which satisfies the equation:
\begin{equation}\label{eqvorticite}
    \partial_t \omega(x,t) + u(x,t)\cdot\nabla \omega(x,t) = 0.
\end{equation}

If the domain $\Omega$ is smooth, then we have a unique global smooth solution of \eqref{Euler}, see \cite{Wolibner1933}. In addition, if $\partial \Omega \in C^{1,1}$ we have the following result due to Yudovitch (see \cite{YUDOVICH19631407}): for any $\omega_0 \in L^1 \cap L^\infty$, there exists a unique solution to \eqref{Euler}, with $u \in L^\infty(\R^+,W^{1,p})$, for every $1<p<\infty$,  and $\omega \in L^\infty(\R^+,L^1 \cap L^\infty)$. The result is true in more general domains, in particular in domains with finite number of corners with angle strictly less than $\pi$ and when the vorticity is compactly supported in $\Omega$, see \cite{LacaveZlatosCorners} and \cite{han2020euler}. In particular, one could take $\Omega$ to be a convex polygon.

Let us denote by $G_\Omega$ the Green's function of $\Omega$. Since the domain is supposed to be simply-connected, the velocity can be recovered from the vorticity through the following Biot-Savart law
\begin{equation}\label{BSG}
u(x,t) =  \int_\Omega \nabla^{\perp}_x G_\Omega(x,y) \omega(y,t) \dd y,
\end{equation}
where $x^\perp = (-x_2,x_1)$.

 The point vortex system is a simplified version of the Euler equations where the vorticity is assumed to be a finite sum of Dirac masses $\omega_0 = \sum_{i=1}^N a_i\delta_{z_i}$. It was introduced by Helmholtz in \cite{Helmholtz1858}, see also \cite{marchioro1993mathematical}. Since  \eqref{eqvorticite} is a transport equation, one expects the vorticity to remain a sum of Dirac masses at some points $z_i(t)$. The  Biot-Savart law \eqref{BSG} reads in this case 
\begin{equation*}
u(x,t) = \sum_{i=1}^N a_i \nabla_x^\perp G_{\Omega}(x,z_i(t)).
\end{equation*}
However, if $x$ is one of the points $z_i(t)$, then this velocity is not defined as $G_\Omega(x,y)$ is singular at $y=x$. But as $x$ approaches $z_i(t)$, the singular part of the velocity defined above is given by fast rotation around that point. More precisely, since the map $G_\Omega - G_{\R^2}$ is harmonic in both its variable on  $\Omega$, the function $\gamma_\Omega : \Omega\times \Omega \rightarrow \R$, $\gamma_\Omega = G_\Omega - G_{\R^2}=G_\Omega-\frac{1}{2\pi}\ln|x-y|$ is smooth. So the singular part of $\nabla_x^\perp G_{\Omega}(x,z_i(t))$ is given by $\frac{(x-z_i)^\perp}{2\pi|x-z_i|^2}$. The point vortex system consists in ignoring this singular part which should have no influence on the motion of $z_i$ itself. Denoting by $\tgm_{\Omega}(x)= \gamma(x,x)$ the Robin function of the domain $\Omega$ we obtain then  the following point vortex dynamic: 
\begin{equation}
\label{ptvortexdynamic}    \forall 1 \le i \le N, \sp \der{z_i(t)}{t} = \sum_{\substack{j= 1 \\ j\neq i}}^N a_j\nabla_x^\perp G_\Omega(z_i(t),z_j(t)) +  a_i\frac{1}{2}\nabla^\perp\tgm_\Omega(z_i(t)) ,
\end{equation}
where the $(a_i)_{1\le i \le N} \in \R\setminus\{0\}$ are the masses of the point vortices $(z_i(t))_{1\le i \le N}$. Equations \eqref{ptvortexdynamic} are called Kirchhoff-Routh equations. 

Due to the singularities of the Green's function, these equations are valid only while the points $z_i(t)$ stay distinct and do not leave the domain $\Omega$.  There exist configurations leading to collapse of the point vortices, but they are exceptional, see \cite{marchioro1993mathematical} for the case of $\R^2$, and \cite{marchioro1984vortex} for the case of the unit disk (we will discuss this more in detail later). 

An important question in fluid dynamics is whether the point vortex system is a good approximation of the Euler equations. There are convergence results in both ways.

Let us first mention that the point vortex system was used as a numerical approximation of the Euler system. More precisely, consider a smooth solution of the Euler equations and construct an initial discrete vorticity which is  a sum of Dirac masses located on a grid $(h j)_{j\in\Z^2}$ where $h\in \R$ is the length of the grid, with masses $h^2\omega_0(h j)$. Solve then the point vortex system with this initial vorticity. In \cite{ConvergencePVtoEulerGoodman}, the authors proved that this point vortex method is consistent, stable, and converges to the smooth solution of the Euler equation. 

The convergence from the Euler system to the point vortex system was also proved in \cite{marchioro1993VorticiesAndLocalization}. More precisely, consider a smooth initial vorticity that is sharply concentrated around some initial point vortices $z_i$: the support of $\omega_0$ is included in the union of the disks of radius $\eps$ around the points $z_i$, with $\eps>0$ being small. The authors of \cite{marchioro1993VorticiesAndLocalization} proved that for any time $\tau$, and for any $\delta >0$, if $\eps=\eps(\tau,\delta)>0$ is small enough then the solution stays sharply concentrated within disks of radius $\delta$ from the points $z_i(t)$ up to time $\tau$. This statement can be seen as a fixed time confinement result.

In this paper we are interested in the so-called \emph{long time confinement problem}, that is we want to know for how long confinement around point vortices remains true. More precisely, we want to understand how $\eps$, $\delta$ and $\tau$ are linked together and we wish to obtain a confinement time $\tau$ as large as possible. We already know from \cite{marchioro1993VorticiesAndLocalization} that $\tau$ goes to infinity when $\eps$ goes to 0, but we would like to obtain an explicit rate as good as possible. 

This problem was already studied by Buttà and  Marchioro \cite{marchiorobutta2018}. These authors assumed that $\delta = \eps^\beta$, with $\beta < 1/2$, and made the following assumptions on the initial vorticity. Assume that $\omega_0\in L^1\cap L^\infty$ and there exists $\nu$ such that
\begin{equation}\label{conditionomega}
    \begin{cases}|\omega_{0}| \le \eps^{-\nu} \\ \ds \omega_0 = \sum_{i=1}^N \omega_{0,i}, \sp  \supp \omega_{0,i} \subset D(z_i,\eps)  \\ 
    \omega_{0,i} \text{ has a definite sign} \\
    \ds \int_\Omega \omega_{0,i} \dd x= a_i.\end{cases}
\end{equation}
Let $\omega(x,t)$ the solution of \eqref{eqvorticite}. 

We denote by $\tau_{\eps,\beta}$ the exit time of the vorticity from the disks of radius $\eps^\beta$:
\begin{equation}\label{defTau}
    \tau_{\eps,\beta} = \sup\left\{t\ge 0, \forall s \in [0,t], \supp \omega(\cdot,s) \subset \bigcup_{i=1}^N D(z_i(s),\eps^\beta) \right\}.
\end{equation} 
For any $N$-tuple of distinct points $(z_i)\in\Omega$, there exists $\eps$ small enough, such that the disks $D(z_i(0),\eps^\beta)$ are disjoints, and therefore this exit time is well defined and strictly positive. The aim is to obtain a lower bound on $\tau_{\eps,\beta}$ depending explicitly on $\eps$. Two results have been obtained in \cite{marchiorobutta2018}. The first is a logarithmic confinement for the whole plane.

\begin{theo}[\cite{marchiorobutta2018}]\label{theologmarchioro} Assume that $\Omega = \R^2$, that the initial vorticity satisfies \eqref{conditionomega} and that the point vortex system with initial data $\sum_{i=1}^n a_i\delta_{z_i}$ has a global solution. Then for every $\beta < 1/2$ there exists $\eps_0 > 0$ and $C >0$ such that
\begin{equation*}
    \forall \eps < \eps_0, \sp \tau_{\eps,\beta} > C |\ln(\eps)|.
\end{equation*}
\end{theo}

The second result is more restrictive, it holds true for the unit disk and for a single point vortex located at the center, but the conclusion is much stronger since it gives a power-law confinement. 
\begin{theo}[\cite{marchiorobutta2018}]\label{theopowermarchioro} Let  $\Omega = D$ and $\omega_0$ satisfying \eqref{conditionomega} with $N = 1$ and $z_1 = 0$, so that it is compactly supported within the disk $D(0,\eps)$. Then for every $\beta < 1/2$ there exists $\eps_0 > 0$ and $\alpha >0$ such that: 
\begin{equation*}
    \forall \eps < \eps_0, \sp \tau_{\eps,\beta} > \eps^{-\alpha}.
\end{equation*}\end{theo}

The aim of this paper is to extend Theorem \ref{theopowermarchioro} to more general domains. We also observe that Theorem \ref{theologmarchioro} can also be extended to bounded domains; we will discuss this problem in a forthcoming paper.

We consider a single point vortex in a simply-connected bounded domain. We assume for simplicity that the mass of the point vortex is 1, but the results below hold true for a general mass. The first question that arises is the location of the point vortex. We will show that there are special points that allow us to obtain the power-law lower bound while for others the logarithmic bound is probably optimal. 

The dynamic of a single point vortex of mass $a$ reduces to the following ODE 
\begin{equation*}
    \der{}{t} z(t) = a\frac{1}{2}\nabla^\perp\tgm_\Omega(z(t)).
\end{equation*}
It is obvious from this ODE that a point vortex is stationary if and only if it is a critical point of the Robin function $\tgm$. The Robin function has been studied, see \cite{gustafsson1979motion}, and we know that such critical points always exist in a bounded domain. Let $x_0$ be a critical point of the Robin function which is fixed for the rest of this paper. From the Riemann mapping theorem we know that there exists a biholomorphic map $T$ from $\Omega$ to the unit disk $D$. We can chose $T$ such that it maps $x_0$ to 0: $T(x_0)=0$. We shall see below that $x_0$ is a critical point of the Robin function if and only if $T''(x_0) = 0$ (see Proposition \ref{propTx0}). This condition therefore characterizes the fact that $x_0$ is a stationary point for the point vortex system. We call such points stationary points.

Our main result is the following.
\begin{theo}\label{theopower} Let $\Omega$ be a simply connected bounded domain of $\R^2$ with $C^{1,1}$ boundary. Let $x_0$ be a stationary point such that $T'''(x_0) = 0$ where $T$ is a biholomorphism from $\Omega$ to the unit disk, mapping $x_0$ to 0. Assume that $\omega_0$ satisfies \eqref{conditionomega} with $N = 1$ and $z_1 = x_0$. Then for every $\beta < 1/2$ and for any $\alpha < \min (\beta, 2-4\beta)$, there exists $\eps_0 > 0$ such that 
\begin{equation*}
    \forall \eps < \eps_0, \sp \tau_{\eps,\beta} > \eps^{-\alpha}.
\end{equation*}
\end{theo}

This extends Theorem \ref{theopowermarchioro} to more general bounded domains. Indeed, in the case of the unit disk we can choose  $T(z)=z$ so the hypothesis given above is verified for the center of the disk. Let us also observe that the hypothesis that $x_0$ is stationary induces no restriction on the domain $\Omega$. Indeed, we recall that Gustafsson \cite{gustafsson1979motion} proved that every simply-connected smooth domain has at least a stationary point. However, the hypothesis that  $T'''(x_0) = 0$  is a condition that not all domains satisfy. We will comment on this in the last section. We will see in particular that any domain which is invariant by some rotation of angle $\theta\in(0,\pi)$ around $x_0$ satisfies the condition $T'''(x_0) = 0$.

In order to understand better the significance of the condition  $T'''(x_0) = 0$, one could assume that the vorticity $\omega$ itself is a point vortex. We study in detail this perturbation problem in Section \ref{confinementfordirac}. We will prove there that if $|T'''(x_0)|<2|T'(x_0)|^3$  then $\tau_{\eps,\beta}=\infty$ if $\eps$ is small enough while if $|T'''(x_0)|>2|T'(x_0)|^3$  then $\tau_{\eps,\beta}$ is in general not better than $C|\ln\eps|$, see Theorem \ref{theodirac}. In other words, in this particular case we have long time confinement better than $C|\ln\eps|$ if and only if  $|T'''(x_0)|<2|T'(x_0)|^3$. However, when $\omega$ is smooth we require the stronger assumption $T'''(x_0) = 0$.

The plan of the paper is the following. In Section \ref{prelim} we introduce some notation and discuss some facts about the Green's function and the point vortex system. In section \ref{confinementfordirac} we consider the particular case when $\omega$ is a point vortex itself. In Section \ref{sect4} we prove Theorem \ref{theopower}. The last section contains some final remarks and some examples of domains for which our theorem applies.

\section{Preliminary tools}\label{prelim}
\small
List of notation:
\begin{itemize}
    \item $\Omega$ is a $C^{1,1}$ bounded and simply connected domain of $\R^2$;
    \item $D(x_0,r)$ is the disk of center $x_0$ and of radius $r$ and $D = D(0,1)$;
    \item $u$ is the velocity of the fluid and $p$ its pressure, satisfying equations \eqref{Euler};
    \item $\omega =  \partial_1u_2 - \partial_2u_1$ is the vorticity of the fluid;
    \item $\supp f$ is the support of the function $f$, namely the closure of the set $\{ x \in \Omega, f(x)\neq 0\}$;
    \item $\delta_z$ is the Dirac mass in $z$;
    \item $G_\Omega$ or $G$ is the Green's function of the domain $\Omega$;
    \item $\gamma_\Omega$ or $\gamma$ is the regular part of $G_\Omega$, see relation \eqref{decompGreen};
    \item $\tgm_\Omega(x) = \gamma_\Omega(x,x)$ is the Robin function;
    \item $C, C_1, C_2, \ldots ; K, K_1, K_2,\ldots, L$, are strictly positive constants that may vary from one line to another, when their value is not important to the result;
    \item $a \cdot b$ is the scalar product of vectors in $\R^2$;
    \item $\nabla f$, $D^2 f$ and $\nabla \cdot g$ are respectively the gradient of $f$, its Hessian matrix, and the divergence of $g$.
\end{itemize}

\subsection{Green's Function}

\normalsize

We recall that the Green's function of a domain $\Omega$ is the solution of
\begin{equation*}
    \Delta_xG_\Omega(x,y) = \delta(x-y)
\end{equation*}
vanishing at the boundary, and at infinity if $\Omega$ is unbounded. It is a symmetric function on $\Omega^2$ that satisfies for $x\neq y$
\begin{equation*}
    \Delta_x ( G_\Omega(x,y)  - G_{\R^2}(x,y)) = 0,
\end{equation*}
which means that $ G_\Omega-G_{\R^2}$ is a function, denoted by  $\gamma_\Omega$, which is harmonic in both of its variable.
Therefore, we have: 
\begin{equation}\label{decompGreen}
    G_\Omega(x,y) = \frac{1}{2\pi}\ln|x-y| + \gamma_\Omega(x,y)
\end{equation} 
where $\gamma_\Omega$ is symmetric and smooth.

In the particular case of $\Omega = D$, we have that
\begin{equation}\label{Greendisque}
    G_{D}(x,y) = \frac{\ln|x-y|}{2\pi} -\frac{\ln|x-y^*||y|}{2\pi}
\end{equation}
where $\ds y^* = \frac{y}{|y|^2}$ is the inverse relative to the unit circle. In particular,
\begin{equation}\label{gammaD}
\gamma_D(x,y)=-\frac{\ln|x-y^*||y|}{2\pi}.
\end{equation}

Let $x_0\in\Omega$. From the Riemann Mapping Theorem (see for instance \cite{ahlfors1966complexanalysis}, chapter 6) and recalling that $\Omega$ is simply-connected, we know that there exists a biholomorphic map $T$ from $\Omega$ to the unit disk $D$. The map $T$ is unique up to compositions with the biholomorphisms of the unit disk which are given by  
\begin{equation*}
\phi_{z_0,\lambda}(z) = \lambda \frac{z_0-z}{1-\overline{z_0}z},
\end{equation*}
with $z_0\in D$ and $|\lambda|=1$ arbitrary. Choosing $z_0$ and $\lambda$ conveniently, we can assume without loss of generality that $T(x_0)=0$ and that $T'(x_0)$ is a strictly positive real number. These two conditions insure the uniqueness of the conformal map $T$. Let us observe that in Theorem \ref{theopower} the mapping $T$ is not unique since we only assume that it maps $x_0$ to 0. However, the condition $T'''(x_0)=0$ does not depend on the choice of $T$ (once we assumed that it maps $x_0$ to 0). Indeed, if $T_1$ and $T_2$ are two biholomorphisms from $\Omega$ to $D$ mapping $x_0$ to 0, then $T_1\circ T_2^{-1}$ is a biholomorphism from $D$ to $D$ mapping 0 to 0. So it must be a rotation: there exists some $\lambda$ of modulus 1 such that $T_1\circ T_2^{-1}(z)=\lambda z$. So $T_1=\lambda T_2$ and therefore $T_1'''(x_0)=\lambda T_2'''(x_0)$. Then $T_1'''(x_0)=0$ if and only if $T_2'''(x_0)=0$.

We will assume from now on that $T(x_0)=0$. The assumption that $T'(x_0)>0$ can be made but is not necessary. In the following, $T'(x_0)$ is a complex number.

The properties of the conformal map $T$ imply a precise description of the Green's function of $\Omega$. Indeed, a Green's function composed with a conformal mapping is another Green's function, see for example \cite{ahlfors1966complexanalysis} chapter 6. Therefore the formula \eqref{Greendisque} yields the following proposition.
\begin{prop}\label{GreenTransportée}
Let $T$ be the biholomorphic mapping introduced above. Then
\begin{equation*}
    G_\Omega(x,y) = G_{D}(T(x),T(y)) = \frac{\ln|T(x)-T(y)|}{2\pi}-\frac{\ln|T(x)-T(y)^*||T(y)^*|}{2\pi}.
\end{equation*}
\end{prop}
In the following, we will have to use both characterizations of the map $T$, as a $\C \rightarrow \C$ map, and as a $\R^2 \rightarrow \R^2$ map. In particular, we recall that $T'(x) = \partial_1 T(x)= \partial_1 T_1(x) + i\partial_1 T_2(x)$, and that $\partial_1 T_1 = \partial_2 T_2$ and $\partial_2 T_1 = -\partial_1 T_2$. So  for any map $f \in C^1(\R^2,\R)$, we have that
\begin{equation*}
    \nabla (f\circ T)(x) = \begin{pmatrix} \partial_1 T_1(x) \partial_1 f(T(x)) + \partial_1 T_2(x) \partial_2 f(T(x)) \\ - \partial_1 T_2(x) \partial_1 f(T(x)) + \partial_1 T_1(x) \partial_2 f(T(x))    \end{pmatrix},
\end{equation*}
or as complex numbers:
\begin{align*}
    \partial_1 (f \circ T)(x) + i \partial_2 (f \circ T)(x) &  = \Re(T'(x))\partial_1 f(T(x)) + \Im(T'(x))\partial_2 f(T(x)) \\ & \sp\sp \sp + i[ - \Im(T'(x))\partial_1 f(T(x)) + \Re(T'(x))\partial_2f(T(x))]   \\
    & =  \overline{T'(x)}(\partial_1 f (T(x))+i\partial_2 f(T(x))) .
\end{align*}
Identifying $\nabla f = \partial_1 f + i \partial_2 f$, we have:
\begin{equation}\label{formulacomposition}
    \nabla (f\circ T)(x) = \overline{T'(x)}\nabla f(T(x)) ,
\end{equation}
where the product above must be understood as the product of two complex numbers. We will frequently use this property in the following.

From  Proposition \ref{GreenTransportée} and from relation \eqref{decompGreen} applied for $\Omega$ and for $D$ we get that
\begin{equation*}
    \gamma_\Omega(x,y) = \gamma_D(T(x),T(y)) + \frac{1}{2\pi}\ln \frac{|T(x) - T(y)|}{|x-y|},
\end{equation*}
and letting $y\to x$ we obtain at the level of Robin functions
\begin{equation}\label{transptgm}
        \forall x\in\Omega, \sp \sp \tgm_\Omega(x) = \tgm_D(T(x)) + \frac{1}{2\pi} \ln|T'(x)|.
\end{equation}
Using the explicit expression of $\gamma_D$, see relation \eqref{gammaD}, we can compute the gradient of $\gamma_\Omega$: 
\begin{equation}\label{gammaexplicite}
\begin{aligned}
\nabla_x\gamma_\Omega(x,y)  
&= \frac{\overline{T'}(x)(T(x)-T(y))}{2\pi|T(x)-T(y)|^2} - \frac{\overline{T'}(x)(T(x)-T(y)^*)}{2\pi|T(x)-T(y)^*|^2} - \frac{(x-y)}{2\pi|x-y|^2}\\
&=\frac{\overline{T'(x)}}{2\pi\overline{(T(x)-T(y))}} - \frac{\overline{T'(x)}}{2\pi\overline{(T(x)-T(y)^*)}} - \frac{1}{2\pi\overline{(x-y)}}
\end{aligned}
\end{equation}
where we used relation \eqref{formulacomposition}.

The mapping $\widetilde{\gamma}_\Omega : \Omega \rightarrow \R^2$, that will be named only $\widetilde{\gamma}$ when there is no ambiguity, has been studied extensively in \cite{gustafsson1979motion}. In particular, it was shown in that paper that  $-\tgm_\Omega$ is a super-harmonic function, that is that $\Delta\tgm_\Omega > 0$, and that it goes to infinity near the boundary like $-\frac{1}{2\pi}\ln(d(x,\partial\Omega))$. This implies  that we have the following proposition (see \cite{gustafsson1979motion}).
\begin{prop}[\cite{gustafsson1979motion}]
For every bounded and simply connected open set $\Omega$, there exists at least one point $x_0$ where $\widetilde{\gamma}$ reaches its minimum.
\end{prop}

Critical points of $\tgm_\Omega$ will be of special interest in the following. The proposition above implies the existence of a critical point of the Robin function. Moreover, if $\Omega$ is convex, then the critical point is also unique (see \cite{gustafsson1979motion}). Though we will not use this result here, it may be interesting to keep this in mind, especially when looking for explicit examples.

\subsection{Point vortex system}

The point vortex dynamic, which is described by equations \eqref{ptvortexdynamic}, can exhibit finite time blow-up of solutions.
One scenario of blow-up is when  two point vortices hit each other in finite time, meaning that there exists $i\neq j$ and $t<\infty$ such that $z_i(t)= z_j(t)$. This phenomena can happen, see for instance \cite{marchioro1993mathematical} or \cite{IftimieMarchioroSelfSimilar} for an example of finite-time collapse of a self similar evolution of point vortices. Another scenario for blow-up is when a point vortex hits the boundary. However, the finite time blow-up is exceptional in the case of  the whole plane, see \cite{marchioro1993mathematical}, and for the unit disk, see \cite{marchioro1984vortex}. In those cases, the $N$-dimensional Lebesgue measure of the set of initial positions that lead to a collapse is 0. The case of a more general bounded domain is an ongoing work.

 Let us recall the convergence theorem obtained in \cite{marchioro1993VorticiesAndLocalization}, for configurations of point vortices that do not lead to finite-time blow-up:
\begin{theo}[\cite{marchioro1993VorticiesAndLocalization}] 
	Let $(a_i, z_i(t))$ be a global solution of the point vortex system \eqref{ptvortexdynamic}.
For every $\delta > 0$ and for every time $\tau$, there exists $\eps > 0$ such that if $\omega_0$ satisfies \eqref{conditionomega} for some $0<\nu < 8/3$, then the vorticity stays confined up to  the time $\tau$ in disks of radius $\delta$ centered on $z_i(t)$.
\end{theo}
In addition, the authors of \cite{Cao2018EulerEO} proved that $\tau_{\eps,\beta} \rightarrow +\infty$ as $\eps$ goes to 0 for any $\beta < 1/3$ (recall that $\tau_{\eps,\beta}$ was defined in relation \eqref{defTau}). However, these theorems don't say anything about how $\delta$ depends on the time $\tau$, or conversely, for how long this confinement remains true, depending on $\eps$. 

For some explicit examples of point vortices we refer to \cite{iftimielargetime}.

As mentioned in the introduction, the stationary point vortices are the critical points of the Robin function. Indeed, the relation \eqref{ptvortexdynamic} with $N=1$ reduces to
\begin{equation*}
    z'(t) = \frac{a}{2}\nabla^\perp\tgm_\Omega(z(t)).
\end{equation*}
We refer to such a critical point $x_0$ by saying that it is a \emph{stationary} point vortex.

Those points can also be characterized in terms of the conformal mapping $T$ as zeroes of $T''$.
\begin{prop}\label{propTx0}
The following conditions are equivalent:
	\begin{itemize}
		\item[(i)] A single point vortex placed in $x_0$ is stationary,
		\item[(ii)]  $\nabla \tgm_\Omega(x_0) = 0$,
		\item[(iii)]  $T''(x_0)=0$.
	\end{itemize}
\end{prop}

\begin{proof}
This proposition was already proved in \cite{gustafsson1979motion}. We 
recall the proof for the convenience of the reader.

We already observed that  $(i)$ and $(ii)$ are equivalent. We only need to prove that $(ii)$ is equivalent to $(iii)$.
	From relations \eqref{transptgm} and \eqref{formulacomposition}, and recalling that $T'$ can't vanish in $\Omega$, we deduce that
	\begin{equation*}
	\nabla \tgm_\Omega(x) = \nabla \tgm_D(T(x)) \overline{T'(x)} + \frac{1}{2\pi}\frac{T'(x) \overline{ T''(x)}}{|T'(x)|^2}.
	\end{equation*}
One can easily check that $0$ is a stationary point for the unit disk, so $\nabla \tgm_D(0)=0$.
	Thus
	\begin{equation*}
	\nabla \tgm_\Omega(x_0) = \frac{1}{2\pi}\overline{\left(\frac{T''(x_0)}{T'(x_0)}\right)}.
	\end{equation*}
Clearly $ \nabla\tgm_\Omega(x_0) = 0$ is equivalent to $T''(x_0) = 0$. 

\end{proof}

\section{Confinement for Dirac mass around a stationary vortex}\label{confinementfordirac}
The aim of this section is to study the case where the vorticity $\omega$  itself is a Dirac mass. The reason we consider this simpler case is because it is easier to have a complete description of what happens. And this in turn gives us an indication of what to expect in the smooth case considered in Theorem \ref{theopower}.

We consider a single point vortex located at $z(t)$ which is close to the stationary point vortex $x_0$. Rescaling time if needed, we can assume without loss of generality that the mass of the point-vortex is 1. The question of knowing if $z(t)$ remains close to $x_0$ is closely related to the notion of  \emph{stability} of $x_0$ that we discuss in what follows.

\subsection{Stability of a stationary point vortex}
Since the mass of the point vortex $z(t)$ is 1, its equation of motion is the following
\begin{equation}\label{singleptvortexdynamic}
z'(t) = \frac{1}{2}\nabla^\perp\tgm_\Omega(z(t)).
\end{equation}
We have that
\begin{equation*}
    \der{}{t}\tgm_\Omega(z(t)) = z'(t)\cdot \nabla \tgm_\Omega(z(t)) = 0
\end{equation*}
which means that the point vortex is evolving on the level set $\widetilde{\gamma}(x)=\mathrm{cst}$. We know from \cite{gustafsson1979motion} that the Robin function $\widetilde\gamma(x)$ goes to infinity as $x$ approaches the boundary of $\Omega$. Therefore, the point vortex $z(t)$ can never reach the boundary so \eqref{singleptvortexdynamic} has a global solution.

Since $-\tgm_\Omega$ is super-harmonic, the eigenvalues of the real symmetric matrix $D^2\tgm_\Omega(x_0)$ have positive sum, meaning that at least one of them is positive. So two main cases can occur: either both eigenvalues are positive, either one is positive and one is negative. We skip the study of the degenerate case when one eigenvalue is 0.

So let $\lambda_+ > 0$ and $\lambda_-$ be the two eigenvalues of $D^2\tgm_\Omega(x_0)$. If $\lambda_- >0 $, the Morse Lemma implies that there exists a change of variables $y=\phi(x)$ in the neighborhood of $x_0$ such that in this neighborhood of $x_0$, 
\begin{equation*}
\tgm_\Omega(x)=    \tgm_\Omega(\phi^{-1}(y)) = \tgm_\Omega(x_0) + (y-x_0)_1^2 + (y-x_0)_2^2.
\end{equation*} 
In particular, in these new local coordinates, the level sets of $\tgm_\Omega$ are circles, so in the real coordinates, they are diffeomorphic to circles. More precisely, the level set $\tgm_\Omega(\phi^{-1}(y)) = \tgm_\Omega(x_0) + r$ is a circle of radius $\sqrt{r}$, provided that $r>0$ is small enough. Because $\phi$ is a $C^\infty$ diffeomorphism, and $\phi(x_0) = x_0$,  there exist constants $k,K>0$ such that $k|y-x_0| < |x-x_0| < K|y-x_0|$ and thus the level set $\tgm_\Omega(x) =\tgm_\Omega(x_0) + r$ is contained in the annulus $\{k\sqrt{r} < |x-x_0| < K \sqrt{r} \}$. 

We conclude from these observations that if $|z(0)-x_0|\le \eps $ with $\eps$ small enough, then we have that $|z(t) - x_0| \le \eps \frac{K}{k}$   for every time $t\ge 0$. This a stability property: if $\eps$ is small enough and if we assume that the point vortex starts at distance $\eps$ of $x_0$, then it remains at distance of order $\eps$ for any time.

Assume now that $\lambda_- <0$. We have in this case that
\begin{equation*}
    \tgm_\Omega(x)=    \tgm_\Omega(\phi^{-1}(y)) = \tgm_\Omega(x_0) - (y-x_0)_1^2 + (y-x_0)_2^2,
\end{equation*} 
so the level sets in a neighborhood of $x_0$ are in the local coordinates $y$ hyperbolas, with the special level set $\tgm_\Omega(x) = \tgm_\Omega(x_0)$ being a union of two line segments. In particular, we will see later that one line segment is repulsive, meaning that when a point vortex evolves on that segment, it moves away from the point $x_0$. We call this situation unstable. In that case, no matter how close the point vortex starts from the critical point $x_0$, it goes away from $x_0$ exponentially fast, as we will see in section \ref{sectionUnstable}.

\begin{figure}[hbt!]
    \centering
    \includegraphics[width = 0.6\textwidth]{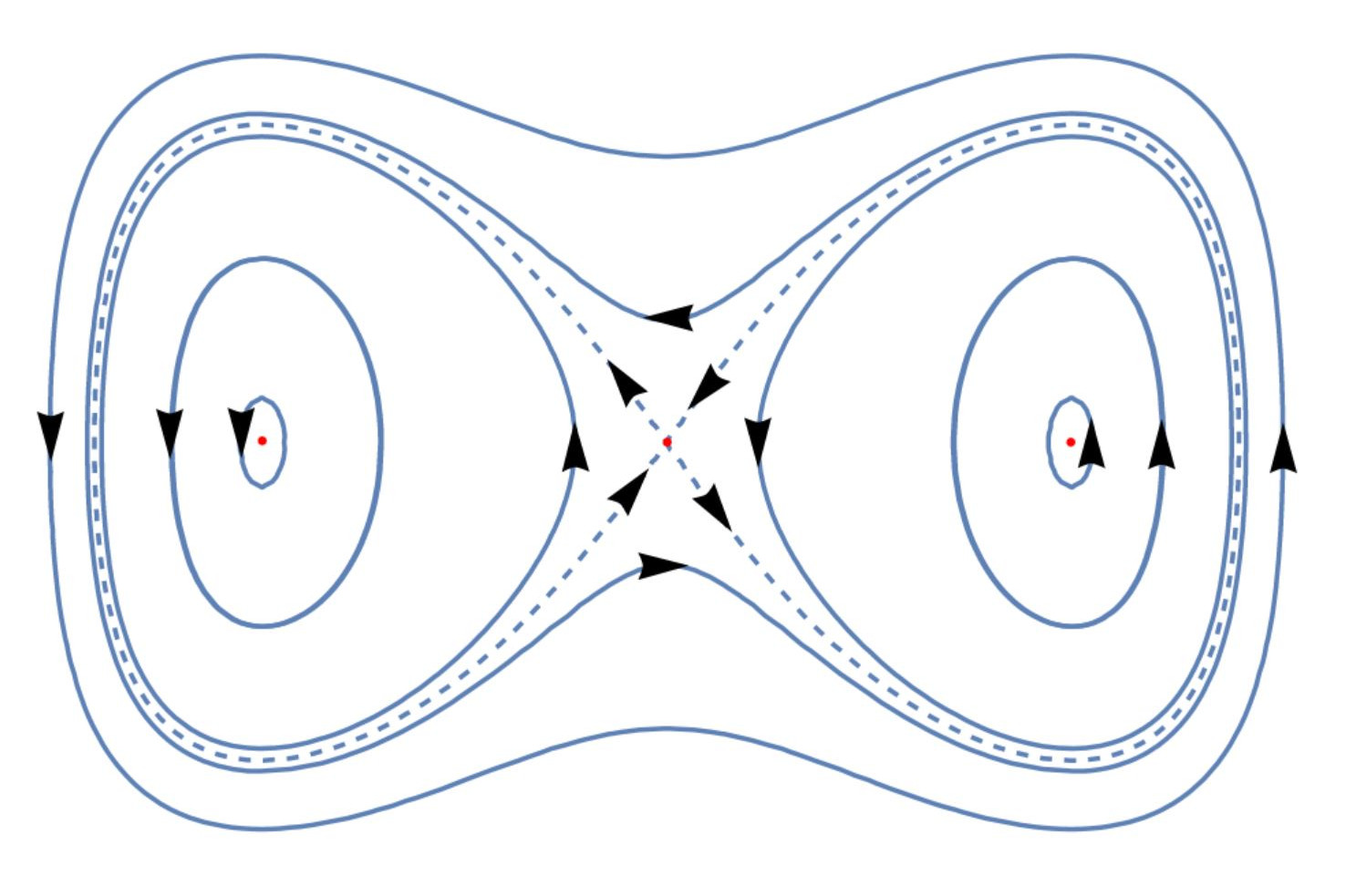}
    \caption{Trajectories of a single point vortex in a domain with two stable points and an unstable one.}
    \label{Traj}
\end{figure}

Let us summarize what we observed above. If the eigenvalues of the matrix $D^2\tgm_\Omega(x_0)$ are both positive, $x_0$ is said to be stable, and a point vortex  close enough to $x_0$ remains indefinitely close to it. If one of the eigenvalues is negative, then $x_0$ is said to be unstable, because there exists a point vortex evolving  on a level set going away from $x_0$. This notion of stability is the same as the one introduced in \cite{marchioro1993mathematical}, Chapter 3. Figure \ref{Traj} shows examples of the two situations. 

\subsection{Taylor expansion near a stationary vortex}

The aim is now to express the stability of a stationary point vortex in terms of the conformal map $T$. We already know from Proposition \ref{propTx0} that $x_0$ is stationary if and only if $T''(x_0)=0$. We will now make a Taylor expansion of $T$ around $x_0$ in order to have a better understanding of the point vortex dynamic near a critical point. We therefore need to compute the Hessian matrix $D^2\tgm_\Omega$ in that point. By relation \eqref{transptgm} we have that
\begin{equation*}
    \frac{1}{2}D^2\tgm_\Omega(x_0)  =  \frac{1}{2}D^2(\tgm_D\circ T)(x_0) + \frac{1}{4\pi}D^2 (\ln|T'|)(x_0) .
\end{equation*}

This expression is actually explicit in terms of the conformal map $T$; indeed, we know that $\tgm_D(x) = -\frac{1}{2\pi}\ln(1-|x|^2)$. We recall that $T''(x_0)=0$, so $\partial_1T'(x_0)=\partial_2T'(x_0)=0$. We recall also that for an holomorphic function $\varphi$ we have that $\partial_1\varphi=\varphi'$ and $\partial_2\varphi=i\varphi'$. Some direct computations give that 
$$
\ds D^2(\tgm_D\circ T)(x_0) = \frac{1}{\pi}\begin{pmatrix}|T'(x_0)|^2 & 0 \\ 0 & |T'(x_0)|^2\end{pmatrix}
$$ 
and 
$$ \partial_i \partial_j (\ln|T'|)(x_0) = \frac{\partial_i \partial_j T'(x_0)\cdot T'(x_0)}
{|T'(x_0)|^2}$$
so that
\begin{equation*}
D^2 (\ln|T'|)(x_0)=\frac1{|T'(x_0)|^2}
\begin{pmatrix}
T'''(x_0)\cdot T'(x_0)&(T'''(x_0))^\perp\cdot T'(x_0)\\
(T'''(x_0))^\perp\cdot T'(x_0)&
-T'''(x_0)\cdot T'(x_0)
\end{pmatrix}.
\end{equation*}
We conclude that
\begin{equation}\label{HessienneTgm}
    \frac{1}{2}D^2\tgm_\Omega(x_0)   = \frac{1}{4\pi}\begin{pmatrix} \ds 2\mu^2 + p & \ds q \\ \ds q & \ds 2\mu^2-p\end{pmatrix},
\end{equation} 
where \begin{equation*} \begin{cases}
    \mu^2 = |T'(x_0)|^2 \vspace{1mm} \\\vspace{1mm}
    \ds p = \frac{1}{\mu^2}T'''(x_0)\cdot T'(x_0) \\
    \ds q = \frac{1}{\mu^2}(T'''(x_0))^\perp\cdot T'(x_0). \end{cases}
\end{equation*}

We compute the characteristic polynomial of the matrix $\frac{1}{2}D^2\tgm_\Omega(x_0)$:
\begin{equation*}
    \det (\frac{1}{2}D^2\tgm_\Omega(x_0)-XI_2)  = X^2 - \frac{\mu^2}{\pi} X + \frac{1}{16\pi^2}\left( 4\mu^4 - p^2- q^2\right)
\end{equation*}
and the eigenvalues are
\begin{equation*}
    \lambda_\pm =  \frac{2\mu^2 \pm \sqrt{p^2 + q^2}}{4\pi}.
\end{equation*}
Furthermore, we notice that $\ds p^2 + q^2 = \frac{|T'''(x_0)|^2}{|T'(x_0)|^2}$. We deduce from the expression of the determinant that $x_0$ is stable when $\lambda_->0$, that is when $2|T'(x_0)|^3 > |T'''(x_0)|$ and unstable when $2|T'(x_0)|^3 < |T'''(x_0)|$.

Another interesting thing to notice is that when $T'''(x_0) = 0$, the matrix $D^2\tgm_\Omega(x_0) = \frac{\mu^2}{2\pi}I_2$, with $I_2$ the identity matrix, and thus the trajectories of the linearized system $z'(t) = (D^2\tgm_\Omega(x_0)(z(t)-x_0))^\perp $ are circles of center $x_0$.

\subsection{Exit time around an unstable stationary point}\label{sectionUnstable}
We now consider the unstable case, meaning $\lambda_+\lambda_- < 0$. We now go back to \eqref{singleptvortexdynamic}. Using relation \eqref{HessienneTgm} and the subsequent calculations, one can check that the eigenvalues of the Jacobian matrix of $\frac{1}{2}\nabla^\perp\tgm_\Omega$ in $x_0$ are $\xi = \sqrt{-\lambda_+\lambda_-}$ and $-\xi$.

We thus know from \cite{hartman1982ordinary} Theorem 6.1 of Chapter 9, and the corollary and remarks associated to this theorem that there exist two local invariant manifolds for the equation, and thus there exist two special solutions to equation \eqref{singleptvortexdynamic} associated with each eigenvalue $\pm \xi$. We denote by $z_1$ the solution that goes away from $x_0$, corresponding to the eigenvalue $\xi$. It satisfies that
\begin{equation}\label{limiteZ2}
    \forall 0<c < \xi, \sp \lim_{t \rightarrow -\infty} |z_1(t) - x_0|e^{-ct} = 0.
\end{equation}

We denote by $t_1$ the first time when $|z_1(t_1)-x_0| = \eps^\beta$ and by $t_2$ the last time when $|z_1(t_2)-x_0| = \eps$. Provided $\eps$ is small enough, we have that $0> t_1 > t_2$. We define $\tilde{z}(t) = z_1(t+t_2)$. This is a solution of equation \eqref{singleptvortexdynamic} such that $|\tilde{z}(0)-x_0| = \eps$ and thus we can set
\begin{equation*}
    \tau_{\eps,\beta} = \max\left\{t\ge 0, \forall s \in [0,t], |\tilde{z}(t)-x_0|\le \eps^\beta \right\}
\end{equation*} which is the exit time as defined in relation \ref{defTau} associated to a point vortex starting at distance $\eps$ of $x_0$ and evolving along the same trajectory as $z_1$. By construction, it satisfies that $\tau_{\eps,\beta}= t_1-t_2 \le -t_2$.

We also know from relation \eqref{limiteZ2} that there exist a constant $M$ and a constant $0<c < \xi$, such that for every $t<0$, we have that $|z_1(t) - x_0|e^{-ct}\le M $. Therefore,
\begin{equation*}
    \eps = |z_1(t_2) - x_0 | \le M e^{ct_2}.
\end{equation*}
Thus, $ \ln \eps \le \ln M + ct_2$ which means that $\tau_{\eps,\beta} \le -t_2 \le \frac{2}{c}|\ln \eps|$ for $\eps$ small enough.

This result empowers the conjecture that one cannot expect any general confinement result, as in Theorem \ref{theologmarchioro}, better than $C |\ln \eps |$, since we have an example of such an exit time in the case where the vorticity is itself a point vortex.

We can summarize the results of this section in the following theorem:
\begin{theo}\label{theodirac}
Assume that $T(x_0) = T''(x_0) = 0$, so that $x_0$ is a stationary point. If  $|T'''(x_0)|<2|T'(x_0)|^3$  then any single point vortex starting close to $x_0$ will remain indefinitely close to it. More precisely, there exists some small constant $\eps_0>0$ such that if $\eps\leq\eps_0$ and $z(0)\in D(x_0,\eps)$ then $z(t)\in D(x_0,C\eps)$ for all times $t$ and $\tau_{\eps,\beta} = +\infty$.

Conversely, the condition $|T'''(x_0)|>2|T'(x_0)|^3$ insures that for every $\eps > 0$, there exists a starting position $z(0)$ such that $|z(0) - x_0| \leq \eps$, but the point vortex  exits the disk $D(x_0,\eps^\beta)$ with $\beta<1$ in finite time $\tau_{\eps,\beta} \le C|\ln\eps|$. 
\end{theo}

\section{Power-law confinement near a stationary vortex}
\label{sect4}

The present section is devoted to the proof of Theorem \ref{theopower}. To simplify the notation, we denote from now on  $\gamma = \gamma_\Omega$, $\tgm = \tgm_\Omega$ and $G = G_\Omega$.

Recall that $\Omega$ is a simply connected and bounded domain, that $T$ is a biholomorphism from $\Omega$ to $D$ such that $T(x_0)=0$. We have that $x_0$ is a stationary point which, according to Proposition \ref{propTx0}, means that $T''(x_0)=0$. 

We assume that $\omega_0$ satisfies \eqref{conditionomega} with $N=1$ and $z_1=x_0$. Without loss of generality we can assume that the mass $a_1$ is 1 (changing the mass means rescaling the time). We therefore have that $\omega_0$ is non negative, is supported in $D(x_0,\eps)$ and $\|\omega\|_{L^1} = 1$.

Let us introduce
\begin{equation}\label{defF}
    F(x,t) = \ds \int \nabla_x^\perp\gamma(x,y)\omega(y,t) \dd y,
\end{equation}
which is the influence of the blob over itself through the boundary. The key point of the proof of Theorem \ref{theopowermarchioro} from \cite{marchiorobutta2018} in the case of the disk is the following property:
\begin{equation}\label{majarticlecercle}
    \forall x,z \in D(0,\delta), \forall t \ge 0, \sp |F(x,t)-F(z,t)| \le C|x-z|\delta^2.
\end{equation}

To prove that the result remains true for more general domains, we will need  a similar inequality. It is important to notice that unlike the case when $\omega$ is a point vortex itself, we don't necessarily have that $F(x_0,t) = 0$ \emph{a priori}. According to \eqref{gammaexplicite} we have that
\begin{equation*}
 2\pi \overline{\nabla_x\gamma(x_0,y)}  = - \frac{T'(x_0)}{T(y)} + \frac{T'(x_0)}{T(y)^*} - \frac{1}{x_0-y}.
\end{equation*}
Recalling that $T(y)^*=1/\overline{T(y)}$ and making the Taylor expansion
\begin{equation*}
T(y)=T'(x_0)(y-x_0)+\frac{T'''(x_0)}6(y-x_0)^3 +\mathcal{O}(|y-x_0|^4)
\end{equation*}
we get that
\begin{align*}
    2\pi \overline{\nabla_x\gamma(x_0,y)}
    & = -\frac{1}{y-x_0}\frac{1}{1 + \frac{T'''(x_0)}{6T'(x_0)}(y-x_0)^2 +\mathcal{O}(|y-x_0|^3)} + T'(x_0)\overline{T}(y) - \frac{1}{x_0-y} \\
    & = -\frac{1}{y-x_0}\left(1 - \frac{T'''(x_0)}{6T'(x_0)}(y-x_0)^2 +\mathcal{O}(|y-x_0|^3)\right) \\ & \hspace{5cm} + T'(x_0)\overline{T'(x_0)(y-x_0)} +\mathcal{O}(|y-x_0|^2) - \frac{1}{x_0-y}.
\end{align*}
We conclude that
\begin{equation}\label{gammax0dev}
    2\pi \nabla_x\gamma(x_0,y) = \overline{\frac{T'''(x_0)}{6T'(x_0)} (y-x_0)} + |T'(x_0)|^2 (y-x_0) + \mathcal{O}(|y-x_0|^2).
\end{equation}

\subsection{Estimate of the influence of the boundary}

The first step of the proof is a technical lemma giving a Lipschitz inequality for $F$. This is the counterpart in $\Omega$ of the relation \eqref{majarticlecercle} valid for the unit disk.

\begin{lemme}\label{lemdev}  We have the following estimate, for $\delta$ sufficiently small and for every $x,y,z \in D(x_0,\delta)$
\begin{equation}\label{dev}
    \overline{\nabla_x\gamma(x,y)} - \overline{\nabla_x\gamma(z,y)} = (x-z)\left(\frac{T'''(x_0)}{6\pi T'(x_0)} + \mathcal{O}\left(\delta\right)\right),
\end{equation}
where the term $\mathcal{O}(\delta)$ is bounded by $C\delta$ with $C$ a constant depending only on $\Omega$ and $x_0$. In particular, if we assume that  $T'''(x_0) = 0$, there exists a constant $K_1=K_1(\Omega,x_0)$ such that
\begin{equation}\label{majdevgamma} \forall x,y,z \in D(x_0,\delta), \sp |\nabla_x\gamma(x,y) - \nabla_x\gamma(z,y)| \le K_1|x-z|\delta. \end{equation}
Furthermore, if we assume that $\supp\omega\subset D(x_0,\delta)$ we have that
\begin{equation}\label{majdev} \forall x,z \in D(x_0,\delta), \sp |F(x,t) - F(z,t)| \le K_1|x-z|\delta. \end{equation}
\end{lemme}

\begin{proof}
Let us define 
$$R(x,y,z) = 2\pi(\overline{ \nabla_x\gamma(x,y) - \nabla_x\gamma(z,y)}).$$
Using  \eqref{gammaexplicite} we can write
\begin{equation*}
        R(x,y,z) = \frac{T'(x)}{T(x)-T(y)} - \frac{T'(x)}{T(x)-T(y)^*} - \frac{1}{x-y} 
        - \left( \frac{T'(z)}{T(z)-T(y)} - \frac{T'(z)}{T(z)-T(y)^*} - \frac{1}{z-y} \right)
\end{equation*}
so
\begin{multline*}
        R(x,y,z) = \frac{T'(x)(T(z)-T(y))-T'(z)(T(x)-T(y))}{(T(x)-T(y))(T(z)-T(y))} + \frac{x-z}{(x-y)(z-y)} \\ - \frac{T'(x)(T(z)-T(y)^*)-T'(z)(T(x)-T(y)^*)}{(T(x)-T(y)^*)(T(z)-T(y)^*)}.
\end{multline*}

We decompose $R$ into two parts:
\begin{equation*}
  R=R_1-R_2
\end{equation*}
with
\begin{equation*}
R_1(x,y,z)=  \frac{T'(x)(T(z)-T(y))-T'(z)(T(x)-T(y))}{(T(x)-T(y))(T(z)-T(y))} + \frac{x-z}{(x-y)(z-y)}
\end{equation*}
and
\begin{align*}
 R_2(x,y,z) &= \frac{T'(x)(T(z)-T(y)^*)-T'(z)(T(x)-T(y)^*)}{(T(x)-T(y)^*)(T(z)-T(y)^*)}\\ 
&= \frac{T'(x)T(z)-T'(z)T(x)}{(T(x)-T(y)^*)(T(z)-T(y)^*)}+ \frac{T(y)^*(T'(z)-T'(x))}{(T(x)-T(y)^*)(T(z)-T(y)^*)}\\
&\equiv R_{21}+R_{22}.
 \end{align*}

Since $T$ is smooth and $T(x_0)=0$, we clearly have that  $|T(x)|\leq C\delta$ and  $|T(y)|\leq C\delta$. So we can estimate
$$
|T(x)-T(y)^*|\geq |T(y)^*| -|T(x)|=\frac1{|T(y)|}-|T(x)|\geq \frac1{C\delta}-C\delta\geq  \frac1{2C\delta}
$$
if $\delta$ is sufficiently small. We can therefore bound $R_{21}$ as follows:
$$
|R_{21}|\leq 4C^2\delta^2 |T'(x)T(z) - T'(z)T(x)|\leq C\delta^2 |x-z|.
$$

Similarly
$$
|T(x)-T(y)^*|\geq |T(y)^*| -|T(x)|=\frac1{|T(y)|}-|T(x)|\geq \frac1{2|T(y)|}
$$
so 
$$
\frac{|T(y)^*|}{|T(x)-T(y)^*|}\leq  2|T(y)||T(y)^*|=2.
$$
Then
$$
|R_{22}|=\frac{|T(y)^*|}{|T(x)-T(y)^*|}\frac{|T'(z) - T'(x)|}{|T(z)-T(y)^*|} 
\leq 4C\delta |T'(z) - T'(x)|\leq C\delta|x-z|^2 \le C \delta^2|x-z|.
$$

We conclude from the estimates above that 
\begin{equation}\label{R2}
 |R_2(x,y,z)|\leq C\delta^2|x-z|.
\end{equation}

To estimate $R_1$, we write it under the form
\begin{equation*}
  R_1(x,y,z)=\frac{N(x,y,z)}{(T(x)-T(y))(T(z)-T(y))(x-y)(z-y)}
\end{equation*}
with
\begin{multline}\label{defN}
    N(x,y,z) = [T'(x)(T(z)-T(y))-T'(z)(T(x)-T(y))](x-y)(z-y) \\ + (x-z) (T(x)-T(y))(T(z)-T(y)).
\end{multline}

As $T$ is holomorphic we observe that $N$ is holomorphic in the variables $x$, $y$ et $z$. One can notice that $N$ is 0 if $x=z$, so it can be factorized by $x-z$. Let us also recall that $\gamma$ is smooth everywhere, so $R$ is smooth too. We proved above that $R_2$ is bounded in $D(x_0,\delta)^3$, so $R_1=R+R_2$ is also bounded in this set. But the denominator of $R_1$ has a factor $(x-y)^2(y-z)^2$, so it follows that $N(x,y,z)$ can be factorized by $(x-y)^2(y-z)^2$. But we observed that $N$ can also be factorized by $x-z$. This implies that there exists a holomorphic function $N_1(x,y,z)$ such that:
\begin{equation}\label{N1}
  N(x,y,z)=(x-y)^2(z-y)^2(x-z)N_1(x,y,z).
\end{equation}
Therefore
\begin{equation}\label{R1div}
\frac{R_1(x,y,z)}{x-z} =\frac{(x-y)(z-y)N_1(x,y,z)}{ (T(x)-T(y))(T(z)-T(y))}. 
\end{equation}

We need now to compute $N_1(x_0,x_0,x_0)$. To do that, we will differentiate five times relation \eqref{N1} and evaluate it in $(x_0,x_0,x_0)$. It is clear that derivatives up to order 4 of $ (x-y)^2(z-y)^2(x-z)$ all vanish at $(x_0,x_0,x_0)$. We can therefore write that
\begin{equation*}
\partial_x^3\partial_z^2 N(x_0,x_0,x_0)= \partial_x^3\partial_z^2 \bigl[ (x-y)^2(z-y)^2(x-z)\bigr](x_0,x_0,x_0)N_1(x_0,x_0,x_0)=  12 N_1(x_0,x_0,x_0).
\end{equation*}

Differentiating relation \eqref{defN} allows to find after some computations that  
\begin{equation*}
    \partial_x^3\partial_z^2 N(x_0,x_0,x_0) = 4T'(x_0)T'''(x_0).
\end{equation*}
Therefore,
\begin{equation*}
N_1(x_0,x_0,x_0)=\frac{ T'(x_0)T'''(x_0)}3
\end{equation*}
so 
\begin{equation*}
N_1(x,y,z)=N_1(x_0,x_0,x_0)+\mathcal{O}(\delta)=  \frac{ T'(x_0)T'''(x_0)}3+\mathcal{O}(\delta).
\end{equation*}

We now go back to \eqref{R1div}. We observe that 
$$
\frac{x-y}{T(x)-T(y)}
$$
is smooth on $D(x_0,\delta)^3$ with value $1/T'(x)$ when $x=y$. Therefore
$$
\frac{x-y}{T(x)-T(y)}=\frac1{T'(x_0)}+\mathcal{O}(\delta).
$$
Similarly
$$
\frac{z-y}{T(z)-T(y)}=\frac1{T'(x_0)}+\mathcal{O}(\delta).
$$
Combining the previous relations results in 
$$
\frac{R_1(x,y,z)}{x-z}=\Bigl(\frac1{T'(x_0)}+\mathcal{O}(\delta)\Bigr)^2\Bigl(\frac{ T'(x_0)T'''(x_0)}3+\mathcal{O}(\delta)\Bigr)
=\frac{ T'''(x_0)}{3T'(x_0)}+\mathcal{O}(\delta)
$$
so 
$$
R_1(x,y,z)= (x-z)\left(\frac{T'''(x_0)}{3T'(x_0)} + \mathcal{O}(\delta)\right).
$$

Recalling \eqref{R2} finally implies that 
\begin{equation*}
    R(x,y,z) = (x-z)\left(\frac{T'''(x_0)}{3T'(x_0)} + \mathcal{O}(\delta)\right),
\end{equation*}
which proves \eqref{dev}. Clearly \eqref{majdevgamma} follows from \eqref{dev} if $T'''(x_0)=0$. Finally, relation  \eqref{majdev} follows from \eqref{majdevgamma} after integrating and recalling that the mass of $\omega$ is 1.
\end{proof}

Comparing \eqref{majdev} and \eqref{majarticlecercle}, we see that we lose the factor $\delta^2$. But in the case $T'''(x_0) = 0$, we still get a factor $\delta$ and this is enough to make our argument work. Actually the proof of \cite{marchiorobutta2018} would still be correct assuming only that $|F(x,t)-F(z,t)| \le C|x-z|\delta$. In our proof, a factor $\delta^2$ would improve the restriction over the power $\alpha$ in Theorem \ref{theopower}. However, please notice that if $T'''(x_0) \neq 0$ we lose the factor $\delta$ and our proof does not work anymore.

\subsection{Estimates of the trajectories}

From now on we assume that $T'''(x_0)=0$.

Let us introduce the center of vorticity: \begin{equation}\label{defbeps}
    B(t) = \int_\Omega x\omega(x,t)\dd x,
\end{equation} and the moment of inertia: \begin{equation*}
     I(t) = \int_\Omega |x-B|^2 \omega (x,t) \dd x.
\end{equation*}

For future needs, let us compute the time derivative of $B$. Recall that $\omega$ satisfies the equation \eqref{eqvorticite} in the sense of distributions and that it is compactly supported. We have that
\begin{align*}
    \der{}{t}{B}(t) & = \int x \partial_t \omega(x,t)\dd x \\
    & = - \int xu(x,t)\cdot \nabla\omega(x,t)\dd x\\
    & = \int (u(x,t)\cdot\nabla)x\ \omega(x,t) \dd x\\
    & = \int u(x,t)\omega(x,t) \dd x\\
    & = \iint\left(\frac{(x-y)^\perp}{2\pi|x-y|^2} + \nabla_x^\perp \gamma(x,y) \right)\omega(y,t)\omega(x,t) \dd x  \dd y
\end{align*} 
where we used  \eqref{BSG} and \eqref{decompGreen}.

Observing that $\frac{(x-y)^\perp}{2\pi|x-y|^2}$ is antisymmetric when exchanging $x$ and $y$ and recalling the definition of $F$ given in  \eqref{defF}, we infer that 

\begin{equation}\label{derB}
 \der{}{t}{B}(t) = \int F(x,t) \omega(x,t) \dd x.
 \end{equation}

Let us define
\begin{equation*}
    R_t = \max\{|x-B(t)|;\ x\in \supp \omega(t) \}
\end{equation*}
and choose some $X(t)\in  \supp \omega(t)$ such that $|X(t)-B(t)|=R_t$. We denote by $s\mapsto X_t(s)$ the trajectory passing through $X(t)$ at time $t$ so that $X_t(t)=X(t)$.

We have the following lemma which allows us to control the time evolution of $R_t$:
\begin{lemme}\label{lemmeVitesseRadiale}
For any $t\leq 	\tau_{\eps,\beta}$ we have that
\begin{equation*}
    \der{}{s}|X_t(s)-B(s)|\big|_{s=t} \le 2 K_1\eps^{\beta} R_t + \frac{5}{\pi R_t^3}I(t) + K_2\left(\eps^{-\nu} \int_{|x-B|>R_t/2}\omega(x,t)\dd x\right)^{1/2}
\end{equation*}
where $\nu$ is the constant from relation \eqref{conditionomega}, $K_1$ is the constant from Lemma \ref{lemdev} and $K_2$ is a universal constant.
\end{lemme}

This lemma shows that in order to obtain upper bounds for the growth of the support of $\omega$, one needs estimates for $I(t)$, $B(t)$ and for the mass of vorticity far from the center of mass. 

\begin{proof}
We have that for any $s\ge 0$ and $t\ge 0$, $X_t'(s) = u(X_t(s),s)$, so
\begin{equation*}
    \der{}{s} |X_t(s)-B(s)| = \left(u(X_t(s),s)-B'(s)\right) \cdot \frac{X_t(s)-B(s)}{|X_t(s)-B(s)|}.
    \end{equation*}
We fix now the time $t \ge 0$, we take $s=t$, and we write $X$  instead of $X_t(t)$. By the Biot-Savart law \eqref{BSG}, the  relation \eqref{defF} and recalling that the vorticity is assumed to be of integral 1, we have that
$$
u(X,t)=F(X,t)+ \int \frac{(X-y)^\perp}{2\pi|X-y|^2}\omega(y,t) \dd y =\int \bigl(F(X,t)+\frac{(X-y)^\perp}{2\pi|X-y|^2} \bigr)\omega(y,t) \dd y .
$$   
Relation   \eqref{derB} now implies that
    \begin{align*}
        \der{}{s} |X_t(s)-B(s)|\big|_{s=t} & = \left[ \int\left( F(X,t)+\frac{(X-y)^\perp}{2\pi|X-y|^2} - F(y,t) \right)\omega(y,t) \dd y \right]\cdot \frac{X-B(t)}{|X-B(t)|}\\
        & = H_1 + H_2
    \end{align*}
where 
\begin{equation*}
    H_1 = \left[ \int\left( F(X,t) - F(y,t) \right)\omega(y,t) \dd y \right]\cdot \frac{X-B(t)}{|X-B(t)|},
\end{equation*}
and
\begin{equation*}
    H_2 = \left[ \int\frac{(X-y)^\perp}{2\pi|X-y|^2}\omega(y,t) \dd y \right]\cdot \frac{X-B(t)}{|X-B(t)|},
\end{equation*}
Thanks to Lemma \ref{lemdev} and recalling that  $t\leq \tau_{\eps,\beta}$, we have that
\begin{align*}
    |H_1| = \left| \int\left( F(X,t) - F(y,t) \right)\omega(y,t) \dd y \right| & \le K_1\eps^\beta \int |X-y|\omega(y,t) \dd y \\
    & \le K_1\eps^\beta2 R_t
\end{align*}
where we used that $\supp\omega\subset \overline{D(B,R)}$. The second term $H_2$ is the same as the left hand side of relation (2.28) in \cite{marchiorobutta2018},  and its estimate is the same: 
\begin{equation*}
|H_{2}| \le\frac{5}{\pi R_t^3}I(t) + \left(\frac{1}{\pi}\eps^{-\nu} \int_{|x-B|>R_t/2}\omega(x,t)\dd x\right)^{1/2}
\end{equation*}
The lemma is now proved.
\end{proof}

\subsection{Estimates of the moments of the vorticity}

We have the following lemma:
\begin{lemme}\label{estimationsmoments} For every $t< \min(\tau_{\eps,\beta},\eps^{-\beta})$, we have:
\begin{equation*}
    I(t) \le K_3\eps^2.
\end{equation*}
and:
\begin{equation*}
    |B(t) - x_0| \le  K_3 \eps,
\end{equation*}
where $K_3$ is a positive constant depending only on $\Omega$ and $x_0$.

\end{lemme}
\begin{proof}
We differentiate $I(t)$:
$$
 I' (t) = \int \left( |x-B|^2 \partial_t \omega (x,t) - 2B'(t)\cdot (x-B)\omega (x,t)\right) \dd x  = \int  |x-B|^2 \partial_t \omega (x,t)  \dd x 
$$
where we used relation \eqref{defbeps}. Next, we use the equation of $\omega$ given in \eqref{eqvorticite} and write
\begin{align*}
    I' (t)
    & =   \int \left( -|x-B|^2 u(x,t)\cdot \nabla\omega (x,t) \right) \dd x \\
    & =   \int  2(x-B) \cdot u(x,t)\omega (x,t) \dd x \\
    & =   \iint 2 (x-B)\cdot \nabla_x^\perp G(x,y) \omega(x,t)\omega(y,t) \dd x\dd y   \\
    & =   \iint 2 (x-B)\cdot \left[ \nabla_x^\perp \gamma(x,y)+ \frac{(x-y)^\perp}{2\pi |x-y|^2} \right] \omega(x,t)\omega(y,t) \dd x\dd y  \text{ } \\
    & =   \iint 2 (x-B)\cdot \nabla_x^\perp \gamma(x,y) \omega(x,t)\omega(y,t) \dd x\dd y  \\
    & \hskip 3cm  + \iint \frac{- x\cdot y^\perp- B\cdot (x - y)^\perp}{\pi |x-y|^2} \omega(x,t)\omega(y,t) \dd x\dd y
\end{align*}    
Exchanging $x$ and $y$ shows that the last term above vanishes. So
\begin{align*}
I' (t)
     & =   \iint 2 (x-B)\cdot \nabla_x^\perp \gamma(x,y) \omega(x,t)\omega(y,t) \dd x\dd y   \\
    & =  \iint 2 (x-B)\cdot [\nabla_x^\perp \gamma(x,y) - \nabla_x^\perp\gamma(B,y)]\omega(x,t)\omega(y,t) \dd x\dd y  \\
    & \leq 2K_1 \eps^{\beta}  \iint  |x-B|^2\omega(x,t)\omega(y,t) \dd x\dd y   \\
    & = 2K_1 \eps^{\beta} I(t)
\end{align*}
where we used \eqref{majdevgamma} with $\delta = \eps^\beta$ because $(x,y,z)\in B(x_0,\eps^\beta)$ (see the definition \eqref{defTau} of $\tau_{\eps,\beta})$. The Gronwall lemma now implies that $\ds I(t) \le I(0) \exp( 2K_1\eps^{\beta} t)$. Since at the initial time we have that $\supp\omega_0$ and $B(0)$ are in the disk $D(x_0,\eps)$ we infer that $I(0) \le 4\eps^2$. We assumed that $t < \eps^{-\beta}$ so we finally obtain that
\begin{equation*}
    \forall t \le \min(\tau_{\eps,\beta},\eps^{-\beta}), \sp \sp I(t) \le 4e^{2K_1}\eps^2.
\end{equation*}

We estimate now the center of vorticity. By relations \eqref{derB} and \eqref{defF}, we have that
\begin{align*}
    \der{}{t}|B(t)-x_0|^2 & = 2 B'(t) \cdot (B(t)-x_0) \\
    & = 2\iint \nabla_x^\perp \gamma(x,y) \omega(y,t) \omega(x,t) \dd x \dd y\cdot (B(t)-x_0). \\
\end{align*}
To estimate this, we use Lemma \ref{lemdev} and relation \eqref{gammax0dev} and we recall that we assume that $T'''(x_0)=0$ and $\supp\omega\in D(x_0,\eps^\beta)$ to obtain
\begin{align*}
    2\pi \nabla_x^\perp \gamma(x,y) & = 2\pi( \nabla_x^\perp \gamma(x,y)- \nabla_x^\perp \gamma(x_0,y)+\nabla_x^\perp \gamma(x_0,y))\\
    & = 2\pi( \nabla_x^\perp \gamma(x,y)- \nabla_x^\perp \gamma(x_0,y)) + |T'(x_0)|^2 (y-x_0)^\perp + \mathcal{O}\left(|y-x_0|^2\right)\\
&=|T'(x_0)|^2 (y-x_0)^\perp +\eps^\beta(|x-x_0|+|y-x_0|)\mathcal{O}(1).
\end{align*}
Putting the estimates above together we have that
\begin{multline*}
\der{}{t}|B(t)-x_0|^2  =  2\iint\left[  \frac{|T'|^2(x_0) (y-x_0)^\perp}{2\pi}+\eps^\beta(|x-x_0|+|y-x_0|)\mathcal{O}(1)\right] \\
\omega(y,t)\omega(x,t) \dd x \dd y \cdot (B(t)-x_0). 
\end{multline*}
But we have the following cancellation:
\begin{align*}
 \iint |T'|^2(x_0) (y-x_0)^\perp \omega(y,t)& \omega(x,t) \dd x\dd y \cdot (B(t)-x_0) \\
&= |T'|^2(x_0)  \iint(y-x_0)^\perp \omega(y,t)\omega(x,t) \dd x \dd y \cdot (B(t)-x_0)\\
&=  |T'|^2(x_0) (B(t)-x_0)^\perp \cdot (B(t)-x_0)\\
&=0.
\end{align*}
Therefore,
\begin{equation*}
\der{}{t}|B(t)-x_0|^2 \le C| B(t)-x_0| \eps^\beta \int |x-x_0| \omega(x,t) \dd x.
\end{equation*}

We notice now that
\begin{align*}
\int |x-x_0|^2\omega(x,t)\dd x & = \int |x-B(t)|^2\omega(x,t)\dd x \\ &\hskip 2cm -\int (x_0 - B(t))\cdot (x-x_0 + x - B(t))\omega(x,t)\dd x \\
& = \int |x-B(t)|^2\omega(x,t)\dd x - (x_0 - B(t))\cdot (B(t)-x_0+B(t)-B(t))  \\
& = I(t)  + |x_0 - B(t)|^2.
\end{align*}

By the Cauchy-Schwarz inequality
\begin{equation*}
\int |x-x_0| \omega(x,t) \dd x\le  \Bigl(\int |x-x_0|^2 \omega(x,t) \dd x \Bigr)^{\frac12}=(I+|B(t)-x_0|^2)^{\frac12},
\end{equation*}
so
\begin{equation*}
\der{}{t}|B(t)-x_0| \le C\eps^\beta (I^{\frac12}+|B(t)-x_0|).
\end{equation*}
By the Gronwall lemma we infer that
\begin{equation*}
    |B(t)-x_0| \le \exp\left( C\eps^\beta t \right)\left(|B(0) - x_0| + C\eps^\beta \int_0^tI(s)^\frac{1}{2} \dd s\right) .
\end{equation*}
We already know that $I(s) \le C\eps^2$, and $|B(0) - x_0| \le \eps$, so 
\begin{equation*}
     |B(t)-x_0| \le \exp(C\eps^\beta t)(\eps +C\eps^\beta t \eps).
\end{equation*}

As we assumed that $\eps^\beta t \le 1$ we conclude that 
\begin{equation*}
     |B(t)-x_0| \le C\eps.
\end{equation*}
This completes the proof of the lemma.
\end{proof}

We have just shown that up to the time  $\min(\tau_{\eps,\beta},\eps^{-\beta})$ the center of mass $B$ stays close to $x_0$  and the moment of inertia $I$ remains small. We need a last technical lemma, which is inspired by the appendix of \cite{Iftimie99Sideris}.

\begin{lemme}\label{est4n}
For every $k \ge 1$ and $t\leq \min(\tau_{\eps,\beta},\eps^{-\beta})$ there exists a small constant $\eps_0=\eps_0(k)$, a large constant $C(k)$ and a constant $K_4$ which depends only on $\Omega$ and $x_0$,  such that if 
\begin{equation*}
\eps\leq\eps_0\qquad\text{and}\qquad    r^4 \ge K_4\eps^{2}\left( 1 + kt\ln(2+t)\right),
\end{equation*}
then
\begin{equation*}
    \int_{|x-B|>r} \omega(x,t)\dd x \le  C(k)\frac{\eps^{k/2}}{r^k}.
\end{equation*}
\end{lemme}

\begin{proof}

Let us introduce the moment of vorticity of order $4n$
\begin{equation*}
    m_n(t) = \int_\Omega |x-B(t)|^{4n}\omega(x,t)\dd x.
\end{equation*}
One can differentiate to obtain, by using relations \eqref{BSG}, \eqref{derB} and recalling that $\nabla \cdot u = 0$
\begin{align*}
    m_n'(t) & = -\int |x-B(t)|^{4n}u(x,t)\cdot \nabla \omega(x,t)\dd x \\
    & \hskip 4cm - 4n\int B'(t)\cdot(x-B(t))|x-B(t)|^{4n-2} \omega(x,t)\dd x  \\
    & =  4n\int u(x,t)\cdot(x-B(t))|x-B(t)|^{4n-2} \omega(x,t)\dd x \\
    & \hskip 4cm - 4nB'(t)\cdot\int (x-B(t)) |x-B(t)|^{4n-2}\omega(x,t)\dd x   \\
    & = 4n\iint \nabla_x^\perp G(x,y)\cdot(x-B(t))|x-B(t)|^{4n-2} \omega(x,t)\omega(y,t)\dd x \dd y \\
    & \sp\sp\sp - 4n\iint \nabla_x^\perp \gamma(z,y) \omega(z,t)\omega(y,t)\dd z \dd y \cdot\int (x-B(t))|x-B(t)|^{4n-2} \omega(x,t)\dd x.
\end{align*}
Recalling that $G(x,y)=\frac{\ln|x-y|}{2\pi}+\gamma(x,y)$ and that $\widetilde{\omega}(x',t) = \omega(x' + B(t),t))$, we can further decompose
\begin{equation*}
 m_n'(t)=a_n(t) + b_n(t) -c_n(t)
\end{equation*}
where
$$ \ds a_n(t) = 4n \iint \frac{(x'-y')^\perp}{2\pi|x'-y'|^2}\cdot x'|x'|^{4n-2} \widetilde{\omega}(x',t)\widetilde{\omega}(y',t)\dd y' \dd x', $$
$$ \ds b_n(t) = 4n \iint \nabla_x^\perp \gamma(x,y)\cdot(x-B(t))|x-B(t)|^{4n-2} \omega(x,t)\omega(y,t)\dd x\dd y,$$
$$ \ds c_n(t) = 4n \iiint \nabla_x^\perp \gamma(z,y)\cdot(x-B(t))|x-B(t)|^{4n-2} \omega(x,t)\omega(y,t) \omega(z,t)\dd x\dd y \dd z.$$

We observe that  $a_n$ is exactly the same quantity that appears in \cite{Iftimie99Sideris} on the last line of page 1726. Observing that the center of mass of  $\widetilde{\omega}(x,t)$ is in 0, we deduce that the estimates given in  \cite{Iftimie99Sideris} are true for $a_n$. More precisely, $a_n$ satisfies the estimate given on the line 5, page 1729 of \cite{Iftimie99Sideris}:
\begin{equation*}
    |a_n(t)| \le C n^2 I(t) m_{n-1}(t).
\end{equation*}
Applying Lemma \ref{estimationsmoments} we obtain that
\begin{equation}\label{estimationAn}
    |a_n(t)| \le C n^2 \eps^2 m_{n-1}(t).
\end{equation}

The terms $b_n$ and $c_n$ are similar. We decompose

\begin{equation*}
    \nabla_x^\perp \gamma(x,y) = \nabla_x^\perp \gamma(x,y)  - \nabla_x^\perp \gamma(x_0,y)+\nabla_x^\perp \gamma(x_0,y).
\end{equation*}
We apply Lemma \ref{lemdev} with $\delta=\eps^\beta$ to deduce that 
$$|\nabla_x^\perp \gamma(x,y)  - \nabla_x^\perp \gamma(x_0,y)|\le K_1|x-x_0|\eps^\beta\leq K_1\eps^{2\beta}$$ 
for all $x,y\in\supp\omega(t)$. Moreover, we also have that $|x-B(t)| \le C\eps^{\beta}$ at least for $\eps$ small enough, since $|x-x_0| \le \eps^\beta$ and $|B(t) - x_0| \le K_3 \eps $. Using these estimates we infer that
\begin{equation*}
    |b_n(t)|\le Cn\eps^{2\beta} \int |x-B(t)|^{4n-1} \omega(x,t)\dd x +d_n\le Cn \eps^{5\beta}m_{n-1}(t) + d_n
\end{equation*}
where
\begin{equation*}
    d_n = 4n\left|\iint \nabla_x^\perp \gamma(x_0,y)\cdot (x-B(t))|x-B(t)|^{4n-2} \omega(x,t)\omega(y,t)\dd x \dd y \right|.
\end{equation*}

Decomposing $\nabla_x^\perp \gamma(z,y) = \nabla_x^\perp \gamma(z,y)  - \nabla_x^\perp \gamma(x_0,y)+\nabla_x^\perp \gamma(x_0,y)$ we obtain that the same estimate holds true for $c_n$:
\begin{equation*}
    |c_n(t)| \le Cn \eps^{5\beta}m_{n-1}(t) + d_n.
\end{equation*}

We estimate now $d_n$. From relation \eqref{gammax0dev} we know that there exists a bounded function $c(y)$ such that $ \nabla_x^\perp \gamma(x_0,y) = [a(y-x_0) + c(y)|y-x_0|^2]^\perp $ where $a=\frac{|T'(x_0)|^2}{2\pi}$. We thus have that
\begin{align*}
    d_n & = 4n\left|\iint [a(y-x_0) + c(y)|y-x_0|^2]^\perp\cdot (x-B(t))|x-B(t)|^{4n-2} \omega(x,t)\omega(y,t)\dd x  \dd y\right| \\
        & \le 4na\left|\int (B(t)-x_0)^\perp\cdot (x-B(t))|x-B(t)|^{4n-2} \omega(x,t)\dd x\right|   \\
        & \hskip 4cm + Cn\iint |y-x_0|^2 |x-B(t)|^{4n-1} \omega(x,t)\omega(y,t) \dd x \dd y\\
    & \le 4na\left|\int (B(t)-x_0)^\perp\cdot (x-B(t))|x-B(t)|^{4n-2} \omega(x,t))\dd x\right|  + C n\eps^{5\beta}  m_{n-1}(t)  \\
    & \le C n (|B(t)-x_0|\eps^{3\beta} + \eps^{5\beta})  m_{n-1}(t).
\end{align*}
Let us recall that Lemma \ref{estimationsmoments} gives that $|B(t)-x_0| \le K_3 \eps$ and therefore
\begin{equation*}
    |b_n(t)| + |c_n(t)| \le Cn \eps^{5\beta}m_{n-1}(t) + 2d_n \le  C n (\eps^{1+3\beta} + \eps^{5\beta})  m_{n-1}(t).
\end{equation*} 
Together with relation \eqref{estimationAn} and recalling that $m_n'(t) = a_n(t) + b_n(t) - c_n(t)$, the last estimate yields that
\begin{equation*}
    |m_n'(t)| \le C n^2 (\eps^2 + \eps^{1+3\beta} + \eps^{5\beta}) m_{n-1}(t).
\end{equation*}

Let us observe now that it suffices to assume that $\beta > 2/5$. Indeed, assume that Theorem \ref{theopower} is proved for any $\beta \in ]2/5,1/2[$. Let $\beta' \le 2/5$ and  $\alpha < \min( \beta', 2-4\beta')$. Since $\alpha<2/5$, one can easily check that there exists $\beta \in ]2/5,1/2[$ such that $\alpha < \min( \beta, 2-4\beta)$ (one can choose $\beta$ close to 2/5). Since $\beta' < \beta$, we also have that $\tau_{\eps,\beta} \le \tau_{\eps,\beta'}$ so $\tau_{\eps,\beta'}>\eps^{-\alpha}$.

We assume in the sequel that $\beta > 2/5$. Due to this additional assumption, we have  that 
\begin{equation}\label{majmnprime}
    |m_n'(t)| \le C n^2 \eps^2 m_{n-1}(t).
\end{equation}
Using Hölder's inequality on $f(x) = \omega^{1/n}(x)$ and $g(x) = |x-B(t)|^{4n-4}\omega^{1-1/n}(x)$ with $p = n$ and $q = \frac{n}{n-1}$, we have that
\begin{align*}
    m_{n-1}(t) & = \int |x-B(t)|^{4n-4}\omega(x,t)\dd x \\
    & \le \left( \int \omega(x,t)\dd x \right)^{1/n}\left( \int|x-B(t)|^{4n}\omega(t,x)\dd x \right)^{(n-1)/n} \\
    & = m_n(t)^{(n-1)/n}.
\end{align*}
This last inequality combined with relation \eqref{majmnprime} gives that
\begin{equation*}
    m_n'(t) \le C n^2 \eps^2 m_n(t)^{(n-1)/n}.
\end{equation*}
We integrate to obtain that
\begin{equation*}
    m_n(t) \le \left( m_n(0)^{1/n} + Cn \eps^2 t\right)^{n}.
\end{equation*}
Clearly, $m_n(0) \le (2\eps)^{4n}$ so, assuming that $\eps\leq1$,
\begin{equation*}
    m_n(t) \le \left( L\eps^2(1 + n t)\right)^{n}
\end{equation*}
for some constant $L$ which depends only on $\Omega$ and $x_0$.
Let us choose any $k\ge 1$, $r$ and $n$ such that
\begin{equation*}
    r^4 \ge 2 L\eps^{2}\left( 1 + k\frac{\ln(2+t)}{\ln2 } t\right)
\end{equation*}
and
\begin{equation*}
    k\frac{\ln(2+t)}{\ln2 } -1 < n \le k\frac{\ln(2+t)}{\ln2 }.
\end{equation*}
This defines $n\ge 1$ since $k\ge 1$ and $\frac{\ln(2+t)}{\ln2 } \ge 1$. It also implies that $2^{n+1} > (2+t)^k$. Thus we have that
{\allowdisplaybreaks
\begin{align*}
    \int_{|x-B|>r} \omega(x,t)\dd x & = \int_{|x-B|>r} \omega(x,t)\frac{|x-B(t)|^{4n}}{|x-B(t)|^{4n}}\dd x \\
    & \le \frac{m_n(t)}{r^{4n}} \\
    & \le \frac{\left( L\eps^2(1 + n t)\right)^{n}}{r^kr^{4n-k}} \\
    & \le \frac{1}{r^k}\frac{\left( L\eps^2(1 + k\frac{\ln(2+t)}{\ln2 } t)\right)^{n}}{\left(2 L\eps^{2}(1 + k\frac{\ln(2+t)}{\ln2 } t)\right)^{n-k/4} }\\
    & = \frac{1}{r^k}\frac{\left( L\eps^2(1 + k\frac{\ln(2+t)}{\ln2 } t)\right)^{k/4}}{2^{n-k/4}}\\
    & = \frac{2^{k/4+1}}{r^k}\eps^{k/2}\frac{\left( L(1 + k\frac{\ln(2+t)}{\ln2 } t)\right)^{k/4}}{ 2^{n+1}}\\
    & \le \frac{2^{k/4+1}}{r^k}\eps^{k/2}\frac{\left( L(1 + k\frac{\ln(2+t)}{\ln2 } t)\right)^{k/4}}{ (2+t)^k}.
\end{align*}
}
The function $\ds t \mapsto \frac{\left( L(1 + k\frac{\ln(2+t)}{\ln2 } t)\right)^{k/4}}{ (2+t)^k}$ is bounded on $\R^+$ for every $k$, so there exists a constant $C(k)$ such that for $\eps$ small enough:
\begin{equation*}
    \int_{|x-B|>r} \omega(x,t)\dd x \le C(k)\frac{\eps^{k/2}}{r^k}.
\end{equation*}
This completes the proof of the lemma.
\end{proof}

\subsection{End of the proof of Theorem \ref{theopower}}

We can now finish the proof of Theorem \ref{theopower}.

We recall that, according to Lemma \ref{lemmeVitesseRadiale}, we have that for each particle such that $|X(t) -B (t)| = R_t$,
\begin{equation*}
    \der{}{s}|X_t(s)-B(s)|\big|_{s=t} \le 2 K_1\eps^{\beta} R_t + \frac{5}{\pi R_t^3}I(t) + K_2\left(\eps^{-\nu} \int_{|x-B|>R_t/2}\omega(x,t)\dd x\right)^{1/2}.
\end{equation*}
Due to the estimates obtained in Lemma \ref{estimationsmoments}, taking $t<\min(\tau_{\eps,\beta},\eps^{-\beta})$, we infer  that
\begin{equation}\label{majorationVitesseRadiale}
    \der{}{s}|X_t(s)-B(s)|\big|_{s=t}  \le  2 K_1\eps^{\beta} R_t + \frac{5K_3\eps^2}{\pi R_t^3} + K_2\left(\eps^{-\nu} \int_{|x-B|>R_t/2}\omega(x,t)\dd x\right)^{1/2}.
\end{equation}
Let us introduce $f$ the solution of the ODE:
\begin{equation}\label{majRtTilde}
    \begin{cases}f'(t) = \ds 4 K_1\eps^{\beta} f(t) + 4\max\left( \frac{5K_3\eps^2}{\pi f^3(t)}, K_2\left(\eps^{-\nu} \int_{|x-B|>f(t)/2}\omega(x,t)\dd x\right)^{1/2} \right) \vspace{1mm}\\
    f(0) = 4\eps. \end{cases}
\end{equation}
We want to show that for every $t \in [0,\min(\tau_{\eps,\beta},\eps^{-\beta})]$, $R_t < f(t) $. Assume that this assertion is false, and let $t_2$ be the first time when it breaks down. Since $f(0) = 4\eps$ and $R_0 \le 2\eps$, we infer that $t_2 > 0$. Let $s \mapsto X_{t_2}(s)$ be a trajectory such that $|X_{t_2}(t_2)-B(t_2)| = R_{t_2}=f(t_2)$. From \eqref{majorationVitesseRadiale} and \eqref{majRtTilde}, we see that
\begin{equation}\label{eqAbsurde}
    \der{}{s}|X_{t_2}(s)-B(s)|\big|_{s=t_2} < f'(t_2).
\end{equation}
However, we have that for every $0< h < t_2 $, 
\begin{equation*}
    |X_{t_2}(t_2-h)-B(t_2-h)| \le R_{t_2-h} < f(t_2-h)
\end{equation*}
which implies that 
\begin{equation*}
    \frac{|X_{t_2}(t_2-h)-B(t_2-h)|-|X_{t_2}(t_2)-B(t_2)|}{-h}  > \frac{f(t_2-h) -f(t_2) }{-h}
\end{equation*}
since $|X_{t_2}(t_2)-B(t_2)| = R_{t_2} = f(t_2)$. Taking the limit as $h \rightarrow 0$, we get a contradiction with \eqref{eqAbsurde}.

We choose now $\alpha$ and $k>6$ such that 
\begin{equation*}
    0<\alpha< \min(\beta,2-4\beta)
\end{equation*}
and 
\begin{equation*}
    k(1/2-\beta) +6\beta - 4 -\nu> 0.
\end{equation*}
We define
\begin{equation*}
    t_2 = \inf\{t>0, f(t) = \eps^\beta\}
\end{equation*}
and
\begin{equation*}
    t_1 = \sup \{ t<t_2, f(t) = \eps^\beta/2 \}
\end{equation*}
so that $t_1<t_2$ and
\begin{equation*}
 f(t_1) = \frac{\eps^\beta}2, \quad  f(t_2) = \eps^\beta\quad \text{and}\quad 
 \frac{\eps^\beta}2 \le  f(t) \le \eps^\beta\quad \forall t \in [t_1,t_2].
\end{equation*}

If $t_2 \geq \eps^{-\alpha}$ then $R_t<f(t)\leq\eps^\beta$ for all $t\in[0,\min(\tau_{\eps,\beta},\eps^{-\beta},\eps^{-\alpha})]=[0,\min(\tau_{\eps,\beta},\eps^{-\alpha})]$. By definition of $\tau_{\eps,\beta}$ we know that $R_{\tau_{\eps,\beta}}=\eps^\beta$. So necessarily $\tau_{\eps,\beta}\geq \eps^{-\alpha}$ which completes the proof of Theorem \ref{theopower}.

We assume from now on that  $t_2 < \eps^{-\alpha}$.

We have the following inequality:
\begin{equation}\label{lemmeEst4nCanBeApplied}
    \forall t \in [t_1,t_2], \sp \sp \left(\frac{f(t)}{2}\right)^4 \ge K_4 \eps^2 (1+kt \ln(2+t))
\end{equation}
which implies that Lemma \ref{est4n} can be applied with $r = f(t)/2$ for $t \in [t_1,t_2]$. Indeed, relation \eqref{lemmeEst4nCanBeApplied} is true since 
\begin{equation*}
    K_4 \eps^2 (1+kt \ln(2+t)) \le K_4 \eps^2 (1+k\eps^{-\alpha} \ln(2+\eps^{-\alpha})) \le \left(\frac{\eps^\beta}{4}\right)^4\le \left(\frac{f(t)}2\right)^4
\end{equation*}
for $\eps$ small enough, as we chose $\alpha < 2 - 4\beta$. Lemma \ref{est4n} yields that
\begin{equation}\label{Lemmaest4nApplied}
    \left(\eps^{-\nu} \int_{|x-B|>f(t)/2}\omega(x,t)\dd x\right)^{1/2}  \le \left( C(k)\frac{\eps^{k/2-\nu}}{(f(t)/2)^k} \right)^{1/2}\qquad\forall t\in[t_1,t_2].
\end{equation}
Since we chose $k$ such that $k(1/2-\beta) + 6\beta - 4-\nu >0$, for $\eps$ small enough we have that
\begin{equation*}
    \eps^{-(k-6)\beta} \eps^{k/2-4-\nu} = \eps^{k(1/2-\beta) + 6\beta - 4-\nu} \le \frac{1}{C(k)2^{2k-6}}.
\end{equation*}
Recalling that $f(t) \ge \eps^\beta/2$ for every $t\in [t_1,t_2]$, we infer that 
\begin{equation*}
    f^{k-6}(t)\ge \Bigl(\frac{\eps^\beta}2\Bigr)^{k-6} \ge C(k) 2^k\eps^{k/2-4-\nu}
\end{equation*}
which in turns gives that
\begin{equation}\label{majParE4F6}
    C(k)\frac{\eps^{k/2-\nu}}{(f(t)/2)^k} \le \frac{\eps^4}{f^6(t)}.
\end{equation}
Using relations \eqref{Lemmaest4nApplied} and \eqref{majParE4F6} in  \eqref{majRtTilde} yields that
\begin{equation*}
    f'(t) \le C \eps^\beta f(t) + C\frac{\eps^2}{f^3(t)}.
\end{equation*}
which implies that
\begin{equation*}
    (f^4)'(t) \le C_1 \eps^\beta f^4(t) + C_1\eps^2
\end{equation*}
for some constant $C_1$.
The end of the argument is now straightforward. We use the Gronwall lemma to obtain that
\begin{equation*}
    f^4(t_2) \le f^4(t_1) e^{C_1\eps^\beta(t_2-t_1)} + \eps^{2-\beta} ( e^{C_1\eps^\beta(t_2-t_1)}-1).
\end{equation*}
Since $t_2-t_1 \le \eps^{-\alpha}$ we have that $C_1\eps^\beta(t_2-t_1)\leq C_1\eps^{\beta-\alpha}<1$ for $\eps$ small enough. Using the inequality  $e^x \le 1 + 2x$ for $0\leq x\leq 1 $ and recalling that $f(t_1) = \eps^\beta/2$ and $f(t_2) = \eps^\beta$ we have that
\begin{equation*}
    \eps^{4\beta} \le (\eps^{\beta}/2)^4 e^{C_1\eps^{\beta-\alpha}} + 2C_1\eps^{2-\alpha}.
\end{equation*}
which implies that
\begin{equation*}
    1 \le \frac{e^{C_1\eps^{\beta-\alpha}}}{16} + 2C_1\eps^{2-4\beta-\alpha} .
\end{equation*}
Since $\alpha < \min(\beta,2-4\beta)$, the right hand side of this inequality goes to $1/16$ as $\eps\to0$. So we obtain a contradiction if $\eps$ is small enough. We thus proved that $t_2 \ge \eps^{-\alpha}$ if $\eps$ is small enough.

In conclusion, for any $\beta<1/2$, and for any $\alpha < \min(\beta, 2- 4\beta)$, there exists $\eps_0$ small enough such that for every $\eps \in (0,\eps_0)$, we have that $\tau_{\eps,\beta} > \eps^{-\alpha}$. This completes the proof of Theorem \ref{theopower}.

\section{Final remarks}

A natural question is the signification of the hypothesis $T'''(x_0) = 0$. As we already discussed, the condition $T''(x_0)=0$ means that $x_0$ is a critical point of the map $\tgm$. Such critical points always exist in a bounded simply connected domain $\Omega$. But to have a strong result of confinement as we proved, we need more than just a critical point of $\tgm$. Indeed, we observed in Section \ref{confinementfordirac} that we should not expect a confinement  time better than $|\ln\eps|$ around unstable points. So we need the stationary point $x_0$ to be at least stable. But  the hypothesis $T'''(x_0) = 0$ is stronger than the stability. Indeed, the stability is characterized by the condition  $|T'''(x_0)| < 2|T'(x_0)|^3 $ which is significantly weaker than $T'''(x_0) = 0$.

So the hypothesis $T'''(x_0) = 0$ is more than just stability. Because we obtained the explicit value of $D^2\tgm(x_0)$, see relation \eqref{HessienneTgm},  we see that this condition is equivalent to the fact $D^2\tgm(x_0)$ is a multiple of the identity. This means that the orbit of a single point vortex in the neighborhood of $x_0$ is almost a circle. We don't know whether this condition is indeed necessary to have a strong confinement as proved in Theorem \ref{theopower}. In our proof we use some crucial cancellations to prove the estimate of the moment of inertia from Lemma \ref{estimationsmoments} that we can't reproduce without the hypothesis $T'''(x_0) = 0$.

One can wonder about the existence of domains satisfying the condition $T'''(x_0) = 0$ in a stationary point $x_0$. We call these domains valid domains. Let us notice that if $T(x_0) = 0$ and $T''(x_0) = 0$, then the condition $T'''(x_0) = 0$ is equivalent to $ (T^{-1})'''(0) = 0$. This means that the image $\Omega=f(D(0,1))$ by any biholomorphic map of the form  $f(z) = x_0 + a_1z + \sum_{k=4}^\infty a_kz^k$ of the unit disk is a valid domain: there exists a biholomorphic map $T:\Omega\to D(0,1)$ such that $T(x_0) = T''(x_0) = T'''(x_0) = 0$ (one can choose $T=f^{-1}$). However, checking that a map $f$ of the form given above is indeed biholomorphic may not be an easy task. 

Let us observe now that regular convex polygons are valid domains. Indeed, by the Schwarz-Christoffel formula (see \cite{ahlfors1966complexanalysis}), there exists a conformal map $f$ from the unit disk to the $n$ sided regular polygon with vertices at $\omega_n = e^{i2\pi/n}$, and its derivative has the form
\begin{equation*}
   f'(z) = c \prod_{k=1}^n \left( 1 - z\omega_n^{-k} \right)^{1-\frac{2}{n}}.
\end{equation*}
We can therefore obtain an explicit value for the second and third derivatives, and check that they vanish at 0. So regular convex polygons are valid domains. 

On the other hand, an ellipse which is not a circle is not a valid domain. Indeed, we know from \cite{Kanas2006} that a conformal mapping from the disk to an ellipse of foci $\pm 1$ mapping 0 to 0 has a Taylor expansion $f(z) = z + A_3z^3 + \ldots$ near 0 with $A_3 > 0$. 

Another method to obtain valid domains is to analyze the effect of rotational invariance. Assume that the domain $\Omega$ is invariant by rotation of angle $\theta\in(0,2\pi)$ around the point $x_0 = 0$. This means that $\tgm_\Omega(x) = \tgm_\Omega(e^{i\theta}x)$ for every $x \in \Omega$. Differentiating this relation and  using relation \eqref{formulacomposition} implies that
\begin{equation*}
    \nabla \tgm_\Omega(e^{i\theta}x)= \begin{pmatrix}\cos \theta & \sin \theta \\ - \sin \theta & \cos\theta \end{pmatrix}\nabla\tgm_\Omega(x)
\end{equation*}
and thus $\nabla \tgm_\Omega(0) = 0$ implying that $0$ is a stationary point. Similarly, differentiating again yields  that
\begin{align*}
    \partial^2_1\tgm_\Omega(x)& =  \partial^2_1[\tgm_\Omega(e^{i\theta x})] = \cos^2 \theta \partial^2_1\tgm_\Omega(e^{i\theta}x) + \sin^2 \theta \partial^2_2\tgm_\Omega(e^{i\theta}x) + 2 \cos\theta \sin \theta \partial_1\partial_2\tgm_\Omega(e^{i\theta}x)\\
    \partial^2_2\tgm_\Omega(x)& = \partial^2_2[\tgm_\Omega(e^{i\theta}x)] = \sin^2 \theta \partial^2_1\tgm_\Omega(e^{i\theta}x) + \cos^2 \theta \partial^2_2\tgm_\Omega(e^{i\theta}x) - 2 \cos\theta \sin \theta \partial_1\partial_2\tgm_\Omega(e^{i\theta}x) \\
    \partial_1\partial_2\tgm_\Omega(x)& = \partial_1\partial_2[\tgm_\Omega(e^{i\theta}x)] = (\cos^2 \theta - \sin^2\theta) \partial_1\partial_2\tgm_\Omega(e^{i\theta}x) \\
    &\hskip 7cm +  \cos\theta \sin \theta ( \partial^2_2\tgm_\Omega(e^{i\theta}x)-\partial^2_1\tgm_\Omega(e^{i\theta}x)). \\
\end{align*}
We set $x=0$ and we subtract the first equation above from the second equation. We get 
\begin{equation*}
    (1-\cos^2\theta + \sin^2\theta)(\partial^2_2\tgm_\Omega(0) - \partial^2_1\tgm_\Omega(0)) = - 4 \cos\theta \sin\theta  \partial_1\partial_2\tgm_\Omega(0).
\end{equation*}
Using this in the third equation yields after some calculations that 
\begin{equation*}
    \theta = \pi \sp\sp\sp \text{ or } \sp\sp\sp \partial_1\partial_2\tgm_\Omega(0) = 0.
\end{equation*}
If $\partial_1\partial_2\tgm_\Omega(0) = 0$ and $\theta\neq\pi$  we observe that  $\partial^2_1\tgm_\Omega(0) = \partial^2_2\tgm_\Omega(0)$. Therefore we have that
$$\theta = \pi \sp \sp \sp \text{ or } \sp \sp \sp \exists \lambda, \sp   D^2\tgm_\Omega(0) = \lambda I_2.$$ 
From \eqref{HessienneTgm} we know that
$$D^2\tgm_\Omega(0) = \frac{\mu^2}{\pi}I_2 + \frac{1}{2\pi}\begin{pmatrix}p & q \\ q & -p\end{pmatrix}.$$ 
Clearly $D^2\tgm_\Omega(0)$ is a multiple of the identity if an only if $p=q=0$. Therefore we have that either $\theta=\pi$ or $p = q = 0$. From the definition of $p$ and $q$ given after relation \eqref{HessienneTgm} we see that for any conformal map $T$ mapping 0 to 0 the condition $p = q = 0$ is equivalent to $T'''(0)=0$. So if $\theta \neq\pi$, then $T'''(0) = 0$ and the domain is therefore valid. This is another proof of the fact that regular polygons are valid domains, and it also gives us many other valid domains.

\begin{figure}[hbt!]
    \centering
    \includegraphics[width = 0.38\textwidth]{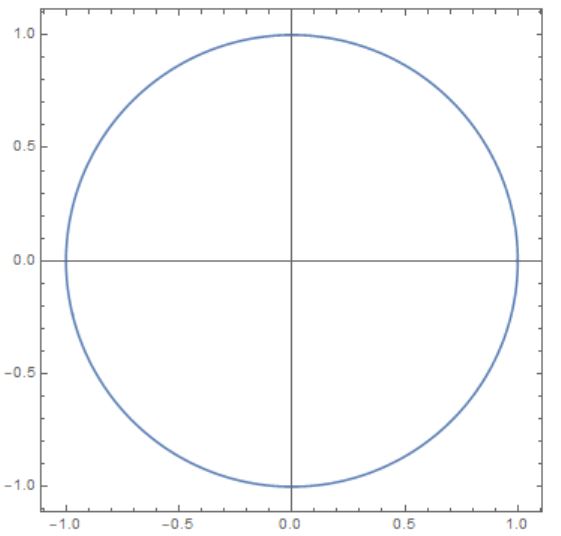}
    \includegraphics[width = 0.38\textwidth]{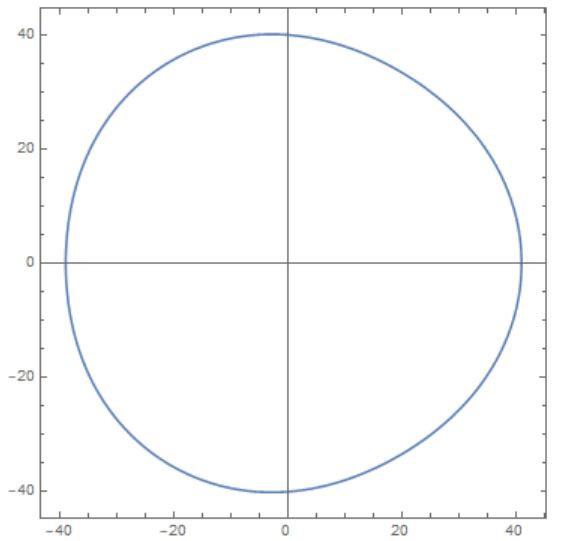}
    \caption{Plot for $f(z) = z$ (left), and $f(z) = 40z + z^4$ (right).}
    \label{figCercle}
\end{figure}

\begin{figure}[hbt!]
    \centering
    \includegraphics[width = 0.38\textwidth]{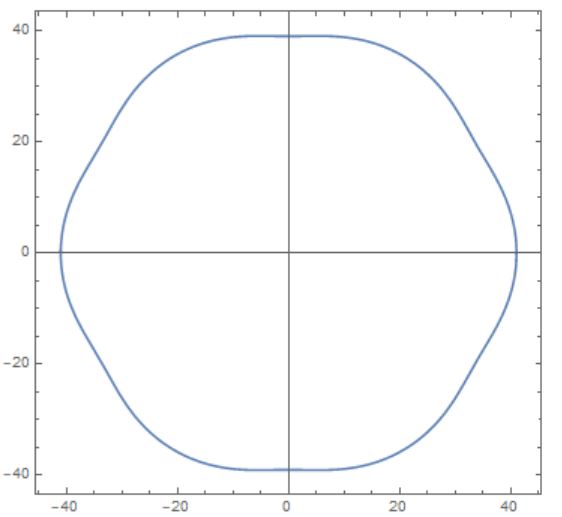}
    \includegraphics[width = 0.38\textwidth]{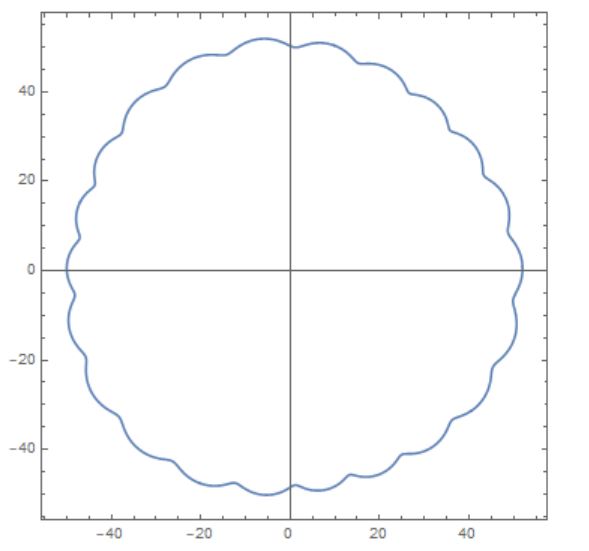}
    \caption{Plot for $f(z) = 40z + z^7$ (left), and $f(z) = 50z + (i+1)z^4+z^{23}$ (right).}
    \label{figModCercle}
\end{figure}

\begin{figure}[hbt!]
    \centering
    \includegraphics[width = 0.38\textwidth]{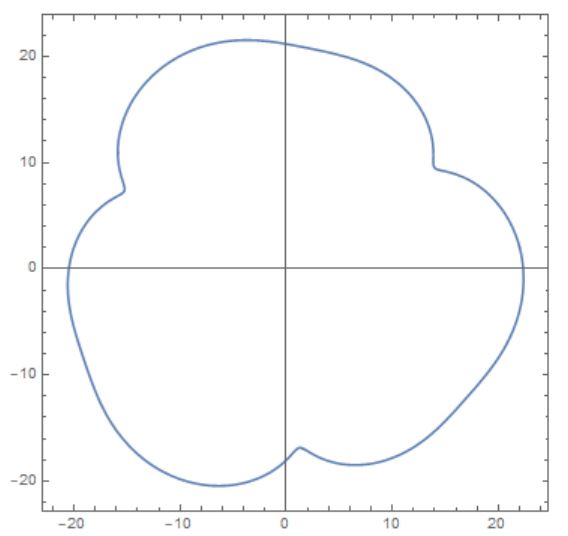}
    \includegraphics[width = 0.38\textwidth]{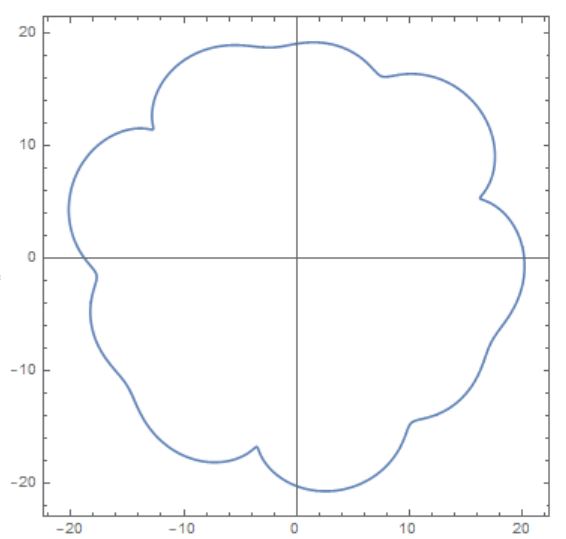}
    \caption{Plot for $f(z) = 20z + (2i+1)z^4+z^7$ (left), and $f(z) = 19z + iz^7 + z^{10}$ (right).}
    \label{figRough}
\end{figure}

\begin{figure}[!h]
    \centering
    \includegraphics[width = 0.38\textwidth]{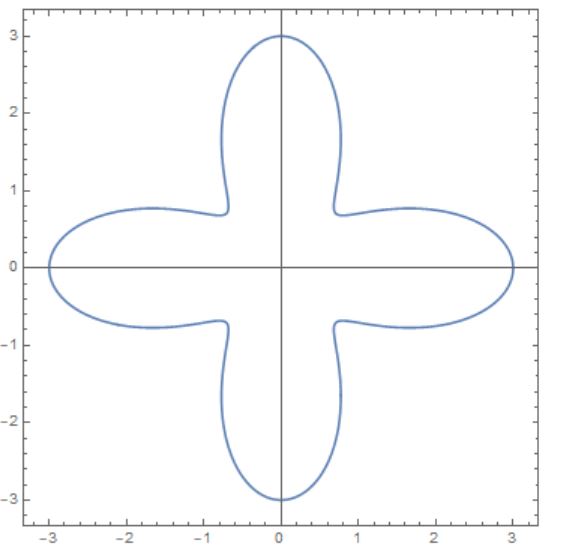}
    \includegraphics[width = 0.38\textwidth]{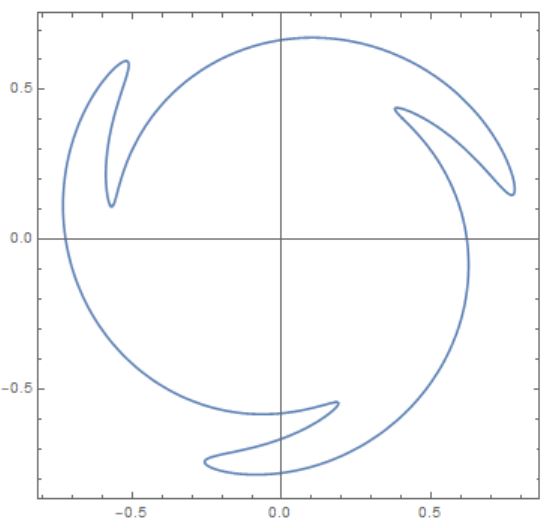}
    \caption{Plot for $b(x) = (2+\cos(8\pi x))e^{i2\pi x}$ (left), and $\ds b(x) = \frac{\exp\left(2i\pi(x+\frac{1}{8}\cos(6\pi x)\right)}{\frac{3}{2}+\frac{1}{4}\cos(6\pi x)}$ (right), with $x \in [0,1]$.}
    \label{figBound}
\end{figure}

We used Mathematica to plot the boundary of some valid domains. We consider different functions $ f(z) = a_1z + \sum_{k=4}^N a_k z^k$, with $N \ge 4$ and $ |a_1| > \sum_{k=4}^N k|a_k|$ in order to ensure that $f$ is injective on the unit disk. We obtain a large class of valid domains, with the unit disk of course and small variations of it (figure \ref{figCercle}), but also larger perturbation of the disk (figure \ref{figModCercle}), and even some very erratic domains (figure \ref{figRough}). Notice that these domains don't necessarily have symmetry properties. Finally, by using the rotational invariance property, we can plot more complicated boundaries without knowing the biholomorphic map, like in figure \ref{figBound}. We plotted images of the interval $[0,1]$ by maps of the form $b(x) = r(x) e^{i2\pi(x+\theta(x))}$, with $r(x) > 0$. We choose $\theta$ and $r$ to be  $1/p$-periodic functions with $p>2$ an integer. This way, $b$ plots a closed curve in $\C$ that is invariant by rotation of angle $2\pi / p \in (0,\pi)$. If this curve does not self intersect and is smooth, then its interior is a valid domain.

\FloatBarrier


\FloatBarrier

\ 

\paragraph{Acknowledgments.}
Dragoș Iftimie has been partially funded by the LABEX MILYON (ANR-10-LABX-0070) of Université de Lyon, within the program “Investissements d’Avenir” (ANR-11-IDEX-0007) operated by the French National Research Agency (ANR). Martin Donati thanks the Dipartimento di Matematica, SAPIENZA Università di Roma for its hospitality. The authors wish to acknowledge helpful discussions with Paolo Buttà and Carlo Marchioro.

\bigskip

\noindent \textbf{Martin Donati}: Université de Lyon, CNRS, Université Lyon 1, Institut Camille Jordan, 43 bd. du 11 novembre,
Villeurbanne Cedex F-69622, France.
\\
Email: donati@math.univ-lyon1.fr

\bigskip

\noindent \textbf{Dragoș Iftimie}: Université de Lyon, CNRS, Université Lyon 1, Institut Camille Jordan, 43 bd. du 11 novembre,
Villeurbanne Cedex F-69622, France.
\\
Email: iftimie@math.univ-lyon1.fr
\\
Web page: http://math.univ-lyon1.fr/\~{}iftimie


\begin{thebibliography}{10}

\bibitem{ahlfors1966complexanalysis}
L.~Ahlfors.
\newblock {\em Complex Analysis}.
\newblock McGraw Hill Publishing Co., New York., 1966.

\bibitem{marchiorobutta2018}
P.~Buttà and C.~Marchioro.
\newblock Long time evolution of concentrated {E}uler flows with planar
  symmetry.
\newblock {\em SIAM J. Math. Anal.}, 50(1):735--760, 2018.

\bibitem{Cao2018EulerEO}
D.~Cao and G.~Wang.
\newblock Euler evolution of a concentrated vortex in planar bounded domains.
\newblock {\em arXiv:1801.01629 [math.AP]}, 2018.

\bibitem{ConvergencePVtoEulerGoodman}
J.~Goodman, T.~Hou, and J.~Lowengrub.
\newblock Convergence of the point vortex method for 2-{D} {E}uler equations.
\newblock {\em Communications on Pure and Applied Mathematics}, 43(3):415 --
  430, 1990.

\bibitem{gustafsson1979motion}
B.~Gustafsson.
\newblock {\em On the Motion of a Vortex in Two-dimensional Flow of an Ideal
  Fluid in Simply and Multiply Connected Domains}.
\newblock Trita-MAT-1979-7. Royal Institute of Technology, 1979.

\bibitem{han2020euler}
Z.~Han and A.~Zlatos.
\newblock Euler equations on general planar domains.
\newblock {\em arXiv:2002.05092 [math.AP]}, 2020.

\bibitem{hartman1982ordinary}
P.~Hartman.
\newblock {\em Ordinary Differential Equations: Second Edition}.
\newblock Classics in Applied Mathematics. Society for Industrial and Applied
  Mathematics (SIAM, 3600 Market Street, Floor 6, Philadelphia, PA 19104),
  1982.

\bibitem{Helmholtz1858}
H.~Helmholtz.
\newblock Über integrale der hydrodynamischen gleichungen, welche den
  wirbelbewegungen entsprechen.
\newblock {\em Journal für die reine und angewandte Mathematik}, 55:25--55,
  1858.

\bibitem{iftimielargetime}
D.~Iftimie.
\newblock Large time behavior in perfect incompressible flows.
\newblock In {\em Partial differential equations and applications}, volume~15
  of {\em S\'{e}min. Congr.}, pages 119--179. Soc. Math. France, Paris, 2007.

\bibitem{IftimieMarchioroSelfSimilar}
D.~Iftimie and C.~Marchioro.
\newblock Self-similar point vortices and confinement of vorticity.
\newblock {\em Communications in Partial Differential Equations},
  43(3):347--363, 2018.

\bibitem{Iftimie99Sideris}
D.~Iftimie, T.~C. Sideris, and P.~Gamblin.
\newblock On the evolution of compactly supported planar vorticity.
\newblock {\em Communications in Partial Differential Equations},
  24(9-10):159--182, 1999.

\bibitem{Kanas2006}
S.~Kanas and T.~Sugawa.
\newblock On conformal representations of the interior of an ellipse.
\newblock {\em Annales Academiae Scientiarum Fennicae. Mathematica},
  31(2):329--348, 2006.

\bibitem{LacaveZlatosCorners}
C.~Lacave and A.~Zlato{\v{s}}.
\newblock {The Euler Equations in planar domains with corners}.
\newblock {\em Archive for Rational Mechanics and Analysis}, 234(1):57--79,
  2019.

\bibitem{marchioro1984vortex}
C.~Marchioro and M.~Pulvirenti.
\newblock {\em Vortex methods in two-dimensional fluid dynamics}.
\newblock Lecture notes in physics. Springer-Verlag, 1984.

\bibitem{marchioro1993mathematical}
C.~Marchioro and M.~Pulvirenti.
\newblock {\em Mathematical Theory of Incompressible Nonviscous Fluids}.
\newblock Applied Mathematical Sciences. Springer New York, 1993.

\bibitem{marchioro1993VorticiesAndLocalization}
C.~Marchioro and M.~Pulvirenti.
\newblock Vortices and localization in euler flows.
\newblock {\em Comm. Math. Phys.}, 154(1):49--61, 1993.

\bibitem{Wolibner1933}
W.~Wolibner.
\newblock Un theor{\`e}me sur l'existence du mouvement plan d'un fluide
  parfait, homog{\`e}ne, incompressible, pendant un temps infiniment long.
\newblock {\em Mathematische Zeitschrift}, 37(1):698--726, 1933.

\bibitem{YUDOVICH19631407}
V.~Yudovich.
\newblock Non-stationary flow of an ideal incompressible liquid.
\newblock {\em USSR Computational Mathematics and Mathematical Physics},
  3(6):1407 -- 1456, 1963.

\end{thebibliography}
\end{document}